# LENSES IN SKEW BROWNIAN FLOW

By Krzysztof Burdzy[1] and Haya Kaspi[2]

*University of Washington and Technion Institute*

We consider a stochastic flow in which individual particles follow skew Brownian motions, with each one of these processes driven by the same Brownian motion. One does not have uniqueness for the solutions of the corresponding stochastic differential equation simultaneously for all real initial conditions. Due to this lack of the simultaneous strong uniqueness for the whole system of stochastic differential equations, the flow contains lenses, that is, pairs of skew Brownian motions which start at the same point, bifurcate, and then coalesce in a finite time. The paper contains qualitative and quantitative (distributional) results on the geometry of the flow and lenses.

**1. Introduction and main results.** The present paper is a continuation of [1] and [4] where an investigation of a stochastic flow of skew Brownian motions driven by a single Brownian motion was initiated. We will study multiple strong solutions to the stochastic differential equation defining the skew Brownian motion. For a fixed starting point, the strong solution to that equation is unique. However, there exist exceptional times ("bifurcation times") when multiple solutions start. We will call pairs of such solutions "lenses" and we will study their properties. Our paper is devoted to a detailed study of a model that belongs to a family of processes analyzed in a series of recent interesting papers by Le Jan and Raimond [9, 10, 11, 12]. We will explain how our model fits into that more general framework at the end of the Introduction.

A skew Brownian motion is a process that satisfies the stochastic differential equation

$$X_t = B_t + \beta L_t, \tag{1.1}$$

Received January 2003; revised October 2003.

[1]Supported in part by NSF Grants DMS-00-71486 and BSF 2000065.

[2]Supported in part by the Fund for the Promotion of Research at the Technion Grant 191-516.

where $B_t$ is a given Brownian motion, $\beta \in [-1, 1]$ is a fixed constant and $L_t$ is the symmetric local time of $X_t$ at 0, that is,

$$(1.2) \qquad L_t = \lim_{\varepsilon \to 0} \frac{1}{2\varepsilon} \int_0^t \mathbf{1}_{(-\varepsilon, \varepsilon)}(X_s) \, ds.$$

The existence and uniqueness of a strong solution to (1.1) and (1.2) was proved by Harrison and Shepp [5]. In the special case of $\beta = 1$, the solution to (1.1) is the reflected Brownian motion. The case $\beta = 0$ is trivial. From now on we will restrict our attention to $|\beta| \in (0, 1)$. An alternative way to define the skew Brownian motion is the following. Consider the case $\beta > 0$. Take a standard Brownian motion $B'_t$ and flip every excursion of $B'_t$ below 0 to the positive side with probability $\beta$, independent of what happens to other excursions. The resulting process has the same distribution as $X_t$ defined by (1.1) (see [6] and [18] for more details).

The following is a straightforward generalization of (1.1). Suppose that $\{B_t, t \in \mathbb{R}\}$ is a Brownian motion on the real line, that is, $\{B_t, t \geq 0\}$ and $\{B_{-t}, t \geq 0\}$ are two independent Brownian motions starting from 0. With probability 1, for all rational $s$ and $x$ simultaneously, the equations

$$(1.3) \qquad X_t^{s,x} = x + B_t - B_s + \beta L_t^{s,x}, \qquad t \geq s,$$

have unique strong solutions, where

$$(1.4) \qquad L_t^{s,x} = \lim_{\varepsilon \to 0} \frac{1}{2\varepsilon} \int_s^t \mathbf{1}_{(-\varepsilon, \varepsilon)}(X_u^{s,x}) \, du.$$

For $s, x, t \in \mathbb{R}$, $t \geq s$, let

$$(1.5) \qquad X_t^{s,x-} = \sup_{\substack{u,y \in \mathbb{Q} \\ u < s \\ X_s^{u,y} < x}} X_t^{u,y},$$

$$(1.6) \qquad L_t^{s,x-} = \sup_{\substack{u,y \in \mathbb{Q} \\ u < s \\ X_s^{u,y} < x}} L_t^{u,y},$$

$$(1.7) \qquad X_t^{s,x+} = \inf_{\substack{u,y \in \mathbb{Q} \\ u < s \\ X_s^{u,y} > x}} X_t^{u,y},$$

$$(1.8) \qquad L_t^{s,x+} = \inf_{\substack{u,y \in \mathbb{Q} \\ u < s \\ X_s^{u,y} > x}} L_t^{u,y}.$$

PROPOSITION 1.1. (i) $X_s^{s,x-} = X_s^{s,x+} = x$ and $X_t^{s,x-} \leq X_t^{s,x+}$ for all $s, x \in \mathbb{R}$ and $t \geq s$, a.s.

(ii) *The processes $t \to X_t^{s,x-}$ and $t \to X_t^{s,x+}$ are Hölder continuous, for all $s, x \in \mathbb{R}$, a.s.*

(iii) *With probability 1, for all $s, x \in \mathbb{R}$ simultaneously, the pairs of processes $(X_\cdot^{s,x-}, L_\cdot^{s,x-})$ and $(X_\cdot^{s,x+}, L_\cdot^{s,x+})$ satisfy* (1.3), *and $L_\cdot^{s,x-}$ and $L_\cdot^{s,x+}$ satisfy* (1.4).

For a fixed "typical" $\omega$ and a "typical" pair $(s, x) \in \mathbb{R}^2$, $X_t^{s,x-} = X_t^{s,x+}$, for all $t \geq s$. This follows easily from the strong uniqueness for a fixed pair $(s, x)$ and the Fubini theorem. Note that the solutions to (1.3) are consistent in the sense that if $X_\cdot^{s,x-} \equiv X_\cdot^{s,x+}$ and $X_\cdot^{u,y-} \equiv X_\cdot^{u,y+}$ for some $u < s$, and $X_s^{u,y-} = x$, then $X_t^{u,y-} = X_t^{s,x-}$ for all $t \geq s$.

This paper is devoted mostly to those $(s, x)$ for which the processes $X_t^{s,x-}$ and $X_t^{s,x+}$ are not identical. We will later show that even if $X_t^{s,x-}$ and $X_t^{s,x+}$ are not identical, there exists some $t_1 = t_1(s, x) < \infty$ such that $X_t^{s,x-} = X_t^{s,x+}$ for $t \geq t_1$.

DEFINITION 1.2. (i) We will say that $\{(X_t^{s,x-}, X_t^{s,x+}), t \in [s, u]\}$, is a *lens* with endpoints $s$ and $u$ if $s < u$, $X_u^{s,x-} = X_u^{s,x+}$ and $X_t^{s,x-} \neq X_t^{s,x+}$ for $t \in (s, u)$. If $\{(X_t^{s,x-}, X_t^{s,x+}), t \in [s, u]\}$, is a lens, then $s$ will be called a *bifurcation time*.

(ii) We will call a bifurcation time $s$ *semi-flat* if for some $s_1 > s$, either $L_t^{s,x-} = L_s^{s,x-} = 0$ for all $t \in [s, s_1]$ or $L_t^{s,x+} = L_s^{s,x+} = 0$ for all $t \in [s, s_1]$. A bifurcation time which is not semi-flat will be called *ordinary*.

It is easy to see that if $\{(X_t^{s,x-}, X_t^{s,x+}), t \in [s, u]\}$, is a lens, then $x = 0$. Hence, we will use the term "bifurcation time $s$" rather than "bifurcation point $(s, x)$." A similar remark applies to the lens endpoint $u$, that is, for every lens $X_u^{s,x-} = X_u^{s,x+} = 0$.

THEOREM 1.3. (i) *With probability 1, the family of all lens endpoints, that is, $u \in \mathbb{R}$ such that for some $s, x \in \mathbb{R}$, $u$ is the endpoint of a lens $\{(X_t^{s,x-}, X_t^{s,x+}), t \in [s, u]\}$, is infinite and countable.*

(ii) *There exist uncountably many bifurcation times, a.s.*

Part (i) of Theorem 1.3 should be clear in view of the following assertion which appears in the next section as Lemma 2.3(i). For any rational times $s_1 < s_2$, with probability 1, the range of $\mathbb{Q} \ni x \to X_{s_2}^{s_1, x}$ consists of two semi-infinite "intervals" $(-\infty, y_1] \cap \mathbb{Q}$ and $[y_2, \infty) \cap \mathbb{Q}$, and a countable set $S \subset (y_1, y_2)$ which does not have accumulation points inside $(y_1, y_2)$.

THEOREM 1.4. (i) *If $|\beta| \in (\frac{1}{3}, 1)$, then semi-flat bifurcation times exist, a.s.*

(ii) *If $|\beta| \in (0, \frac{1}{3})$, then there are no semi-flat bifurcation times, a.s.*

The critical exponent $\frac{1}{3}$ appeared in Corollary 1.5 of [4], which says that there exist random times $t$ when $\beta L_t^{0,0} = \sup_{s \leq t} B_s$ if and only if $\beta > \frac{1}{3}$. It would be interesting to find a direct link between that result and Theorem 1.4 above, for example, via a time reversal argument; so far, we are unable to provide such a direct link.

The next result is concerned with solutions to (1.3). For a deterministic function $t \to B_t$, the pair of equations (1.3) and (1.2) have a clear meaning, not depending on any probabilistic concepts (but that does not mean that a solution must exist for every deterministic $B_t$). Hence, for any "fixed" Brownian paths $B_t$, we can consider all solutions to (1.3) and (1.2) with a given initial condition $(s, x)$.

THEOREM 1.5. (i) *With probability* 1, *there exist* $(s, x) \in \mathbb{R}^2$ *where three distinct solutions to* (1.3) *start.*

(ii) *There are no* $(s, x)$ *where four distinct solutions start.*

We believe that "ordinary" bifurcation times are typical and "semi-flat" bifurcation times are less typical. This informal claim is supported by Theorems 1.3(ii) and 1.4(ii). One can probably formalize the claim by computing the Hausdorff dimensions of ordinary and semi-flat bifurcation times for various values of $\beta$. We will not do this in the present paper. Instead, we will focus on a subfamily of bifurcation times and the corresponding lenses because we can give several fairly explicit formulas in this special case.

We will write $X_t^{x-}$ instead of $X_t^{0,x-}$ and, similarly, $X_t^{x+} = X_t^{0,x+}$. The corresponding local times will be denoted $L_t^{x-}$ and $L_t^{x+}$. Recall that $X_t^{0,x}$, $x \in \mathbb{Q}$, denotes the family of unique strong solutions to (1.3). For rational $x$, we will write $X_t^x = X_t^{0,x}$ and $L_t^x = L_t^{0,x}$. We will call a bifurcation time $s$ *anticipated* if it corresponds to a lens $\{(X_t^{s,0-}, X_t^{s,0+}), t \in [s, u]\}$, and for some $y \in \mathbb{R}$, we have $X_s^{y+} = 0$ and $X_s^{y-} = 0$. In other words, an anticipated bifurcation point may appear only on the trajectory of one of the processes $X_t^{y-}$ or $X_t^{y+}$ for some real $y$. Note that $X_t^{0-} \equiv X_t^{0+}$, a.s., and that for every $x \neq 0$, there exists a random $t_1 > 0$ such that $L_t^{x-} = L_t^{x+} = 0$ for all $t \in [0, t_1]$.

If $s$ is a bifurcation time, $U_s$ will denote the lens $\{(X_t^{s,0-}, X_t^{s,0+}), t \in [s, u]\}$, shifted to 0, that is,

$$U_s = \{(X_{s+t}^{s,0-}, X_{s+t}^{s,0+}), t \in [0, u-s)\}.$$

We let $U_s(t) = \Delta$, a cemetery state, for $t \geq u - s$. Let $\widehat{L}_t^s = L_t^{s,0-} + L_t^{s,0+}$, $\sigma_t = \inf\{u : \widehat{L}_u^s > t\}$, and $Z_t^s = \beta |L_{\sigma_t}^{s,0-} - L_{\sigma_t}^{s,0+}|$. In other words, $Z_t^s$ is the distance between $X_t^{s,0-}$ and $X_t^{s,0+}$ on the time scale defined by the local time clock. Let $\ell^s = \inf\{t > 0 : X_t^{s,0-} = X_t^{s,0+}\}$ and $\ell_Z^s = \inf\{t > 0 : Z_t^s = 0\}$.

Brownian motion is continuous so $L_t^{s,0-}$ and $L_t^{s,0+}$ increase on disjoint intervals for $t \in (s, \ell^s)$, whose endpoints have no accumulation points inside $(s, \ell^s)$. This and the definition of $Z_t^s$ show that on some intervals $Z_t^s$ increases at the rate $\beta$ and on some other intervals it decreases at the rate $\beta$. In other words, it is a piecewise linear function with the slope $\beta$ or $-\beta$ almost everywhere, on the interval $[0, \ell_Z^s]$. Let $J_t^s = 0$, if at time $\sigma_t$, $X_t^{s,0-}$ is at 0, and $J_t^s = 1$, if at time $\sigma_t$, $X_t^{s,0+}$ is at 0. If $\beta > 0$, then $J_t^s$ is the indicator function of the intervals where $Z_t^s$ is increasing.

Let $Q^{x,y}$ denote the distribution of $\{(X_t^{0,x}, X_t^{0,y}), t \geq 0\}$ killed at the time $\zeta = \inf\{t > 0 : X_t^{0,x} = X_t^{0,y}\}$. Note that $\zeta < \infty$, a.s., by the result in [1]. Although we have defined $X_t^{0,x}$ and $X_t^{0,y}$ for rational $x$ and $y$ only, it is clear that the definition of the distribution $Q^{x,y}$ applies to any real $x$ and $y$.

The next theorem involves a $\sigma$-finite measure $Q$ on $C[0, \infty)^2$. We will now introduce some notation related to this measure. The measure $Q$ is supported on pairs of trajectories, say, $(\{X_t, t \geq 0\}, \{\widetilde{X}_t, t \geq 0\})$. There exists a process $\{B_t, t \geq 0\}$ (whose trajectories have Brownian path properties away from $t = 0$), and

$$X_t = B_t + \beta L_t, \qquad \widetilde{X}_t = B_t + \beta \widetilde{L}_t,$$

where

$$L_t = \lim_{\varepsilon \to 0} \frac{1}{2\varepsilon} \int_0^t \mathbf{1}_{(-\varepsilon, \varepsilon)}(X_s)\, ds, \qquad \widetilde{L}_t = \lim_{\varepsilon \to 0} \frac{1}{2\varepsilon} \int_0^t \mathbf{1}_{(-\varepsilon, \varepsilon)}(\widetilde{X}_s)\, ds.$$

We write $\widehat{L}_t = L_t + \widetilde{L}_t$, $\widehat{\sigma}_t = \inf\{s : \widehat{L}_s > t\}$, and $Z_t = \beta|L_{\widehat{\sigma}_t} - \widetilde{L}_{\widehat{\sigma}_t}|$. Hence, $Z_t$ is the distance between the two components of the $Q$-lens on the time scale defined by the local time clock. Finally, we let $\ell = \inf\{t \geq 0 : X_t = \widetilde{X}_t\}$ and $\ell_Z = \inf\{t \geq 0 : Z_t = 0\}$.

THEOREM 1.6. (i) *Let $G$ denote the set of all anticipated bifurcation times. With probability 1, all anticipated bifurcation times are times when the Brownian motion $B_t$ attains its running extremum, that is, if $s \in G$, then $B_s = \sup_{t \leq s} B_t$ or $B_s = \inf_{t \leq s} B_t$. The set $G$ is countable.*

(ii) *There exists a unique (up to a multiplicative constant) $\sigma$-finite measure $Q$ on $C[0, \infty)^2$ which is Markov on every interval $(s, \infty)$, $s > 0$, with the transition probabilities $Q^{x,y}$, and such that both paths start from 0, $Q$-a.e. We have $\lim_{y \downarrow 0}(1/y)Q^{-y,0} = cQ$, for a constant $c$. We will normalize $Q$ so that $c = 1$ in the last formula.*

(iii) *Let $|C|$ denote the Lebesgue measure of $C \subset \mathbb{R}$. For a suitable normalization of $Q$, (nonrandom) Borel sets $A \subset \mathbb{R}$ and bounded continuous functions $f : C[0, \infty)^2 \to \mathbb{R}$,*

$$(1.9) \quad E \sum_{s \in G} \mathbf{1}_A(B_s) f(U_s) = \left( \frac{1-\beta}{1+\beta} |A \cap (-\infty, 0]| + |A \cap (0, \infty)| \right) \int f\, dQ.$$

(iv) *Let $\mathcal{B}$ be the collection of pairs $(B_s, U_s)$ in $\mathbb{R} \times C[0,\infty)^2$ for all $s \in G$. The point process $\mathcal{B}$ is not Poisson.*

(v) *Let $Q^{Z,J}$ be the $Q$-distribution of the process $(Z_t, J_t)$. Assume that $\beta > 0$ and let $\mathcal{A}$ be the collection of all pairs $(B_s, \{(Z_t^s, J_t^s), t \in [0, \ell_Z^s]\})$, where $s \in G$. Let $\mathcal{D}$ be the space of cadlag functions mapping a finite or infinite interval $[0, \zeta]$ to $[0, \infty) \times \{0, 1\}$. Then $\mathcal{A}$ is a Poisson point process on $\mathbb{R} \times \mathcal{D}$ with intensity measure*

$$\mathbf{1}_{(-\infty, 0]}(x)\, dx \times \frac{1-\beta}{1+\beta} Q^{Z,J} + \mathbf{1}_{(0,\infty)}(x)\, dx \times Q^{Z,J}.$$

*An analogous result holds in the case $\beta < 0$, by symmetry.*

(vi) *For any fixed $b > 0$, the $Q^{Z,J}$-distributions of $\{Z_t, t \in [0, \ell_Z]\}$ and $\{Z_{\ell_Z - t}, t \in [0, \ell_Z]\}$, conditional on $\{\sup_{t \geq 0} Z_t = b\}$, are identical.*

(vii) *For any fixed $b > 0$, the $Q$-distributions of processes $\{(X_t, \widetilde{X}_t), t \in [0, \ell]\}$ and $\{(X_{\ell - t}, \widetilde{X}_{\ell - t}), t \in [0, \ell]\}$, conditional on $\{\sup_{t \geq 0} Z_t = b\}$, are different.*

The measure $Q$ is the "distribution" of a lens $U_s$. Theorem 1.6 shows that the point process of lenses has only some of the properties of the familiar excursion processes. It is not a Poisson point process, but one can find a Poisson point structure by restricting attention to a functional of a lens, namely, $(Z_t^s, J_t^s)$. The point process of lenses does satisfy a Maisonneuve-type formula (1.9) (cf. [13]). A similar remark applies to the "lens law" $Q$. The lens law $Q$ is not invariant under time-reversal, but the functional $Z_t$ of a lens is invariant under this transformation.

As a byproduct of the proof of Theorem 1.6, we obtain the following Williams-type decomposition of $Z_t$ under $Q$ [see the remark before the proof of Theorem 1.6(vi) for an alternative presentation]. The process $Z_t$ is piecewise linear on every closed interval contained in $(0, \ell_Z)$. The slope of $Z_t$ is either $\beta$ or $-\beta$, at almost every $t$. Suppose that $\beta > 0$, condition the process $Z_t$ on $\{\sup_{t \geq 0} Z_t = b\}$ and let $\nu$ be such that $Z_\nu = b$. For $t < \nu$, the process $Z_t$ changes the slope from $\beta$ to $-\beta$ at the rate $(1-\beta)/(2Z_t)$, and from $-\beta$ to $\beta$ at the rate $(1+\beta)/(2Z_t)$. By Theorem 1.6(vi), the evolution of $Z_t$ for $t > \nu$, may be described using time-reversal. Moreover, $\{Z_t, t \in [0, \nu]\}$ and $\{Z_{\ell_Z - t}, t \in [\nu, \ell_Z]\}$ are independent under $Q$ given $\{\sup_{t \geq 0} Z_t = b\}$.

The rest of the paper contains some additional results and the proofs of the main theorems; it is divided into two more sections. The next section deals with the definition and properties of the flow and the existence of ordinary and semi-flat bifurcation times. The last section is devoted to anticipated bifurcation times and their distributions.

We will now explain how some of our results can be derived from those of [9, 10, 11, 12], although we will use our own elementary methods in the formal proofs to keep our paper self-contained. The semigroup corresponding to the skew Brownian motion is symmetric. By [9], it is possible to construct a coalescing flow $\varphi_{s,t}$ such that for every function $f$ in the domain of the generator of skew Brownian motion, for all $x$ and $s < t$,

$$f(\varphi_{s,t}(x)) = f(x) + \int_s^t f'(\varphi_{s,u}(x))\,dB_u + \tfrac{1}{2}\int_s^t f''(\varphi_{s,u}(x))\,du.$$

This shows that $\varphi_{s,t}(x) = x + B_t - B_s + \beta L_t^{s,x}$, a property analogous to (1.3). For $s = 0$, any $n \geq 1$, and any $x_1, x_2, \ldots, x_n$, one can solve (1.3) simultaneously for all initial conditions $x_1, x_2, \ldots, x_n$. One can show that this $n$-point motion is Feller and then one can apply a result from [10] to prove that there exists a unique, up to a modification, coalescing flow solving (1.3). Our processes $X_t^{s,x-}$ and $X_t^{s,x+}$ are cadlag and caglad modifications of the flow, in the space variable. Lemma 2.2 can be deduced from the flow property. Lemma 2.5 follows from the fact that the flow is coalescing.

**2. Skew Brownian motion flow.** We fix some $\beta \in (-1,0) \cup (0,1)$ in this section, until stated otherwise.

Recall the notation and definitions from Section 1. The modulus of continuity $\delta_{[a,b]}(r)$ of the Brownian path on the interval $[a,b]$ is defined by

$$\delta_{[a,b]}(r) = \sup\{|B_t - B_s| : s, t \in [a,b], |s - t| \leq r\}.$$

Note that $r \to \delta_{[a,b]}(r)$ is nondecreasing.

PROPOSITION 2.1. *With probability* 1, *for all rational* $s, x$ *with* $s \leq a$, *and all* $u, v \in [a, b]$,

$$|X_u^{s,x} - X_v^{s,x}| \leq 2\delta_{[a,b]}(|u - v|).$$

PROOF. Since rationals are countable, it is enough to prove the proposition for fixed $s$ and $x$. Let $X_t^{s,x}$ and $L_t^{s,x}$ be as in (1.3). For any $t > 0$, let $g(t) = \inf\{u : L_u^{s,x} = L_t^{s,x}\}$ and $d(t) = \sup\{u : L_u^{s,x} = L_t^{s,x}\}$. We will argue that $X_{g(t)}^{s,x} = 0$. Suppose otherwise. Then for some $\varepsilon > 0$, $X_t^{s,x} \neq 0$ for all $u \in (g(t) - \varepsilon, g(t) + \varepsilon)$. It follows from (1.2) that $L_u^{s,x} = L_{g(t)}^{s,x}$ for all $u \in (g(t) - \varepsilon/2, g(t) + \varepsilon/2)$, and this contradicts the definition of $g(t)$. Similarly, $X_{d(t)}^{s,x} = 0$, a.s. We have $X_{g(t)}^{s,x} = 0$ and $X_{d(t)}^{s,x} = 0$ for all rational $t$ simultaneously, a.s. Then it is easy to see that, in fact, $X_{g(t)}^{s,x} = 0$ and $X_{d(t)}^{s,x} = 0$ for all *real* $t \geq 0$ simultaneously, a.s. Consider $u, v \in [a, b]$, and assume without loss of generality that $u < v$. If $d(u) \geq g(v)$, then $|\beta L_u^{s,x} - \beta L_v^{s,x}| = 0$. Otherwise, $|d(u) - g(v)| \leq |u - v|$ and

$$|\beta L_u^{s,x} - \beta L_v^{s,x}| = |\beta L_{d(u)}^{s,x} - \beta L_{g(v)}^{s,x}| = |B_{d(u)} - B_{g(v)}|$$
$$\leq \delta_{[a,b]}(|d(u) - g(v)|) \leq \delta_{[a,b]}(|u - v|).$$

Hence, $\beta L_t^{s,x}$ has the same modulus of continuity as $B_t$ on $[a,b]$, or smaller one. Since $X_\cdot^{s,x}$ is the sum of two functions $x + B_\cdot - B_s$ and $\beta L_\cdot^{s,x}$ with the modulus of continuity bounded by $\delta_{[a,b]}(\cdot)$, its modulus of continuity is bounded by $2\delta_{[a,b]}(\cdot)$. □

Since for every $\alpha < \frac{1}{2}$ the Brownian motion is $\alpha$-Hölder continuous, the same is true for $X_t^{s,x}$, for all rational $s$ and $x$ simultaneously.

LEMMA 2.2. *With probability 1, for all rational $s_1, s_2, x_1, x_2$ simultaneously, if $s_1 \leq s_2$, then either $X_t^{s_1,x_1} \leq X_t^{s_2,x_2}$ for all $t \geq s_2$ or $X_t^{s_1,x_1} \geq X_t^{s_2,x_2}$ for all $t \geq s_2$.*

PROOF. The claim is a part of Proposition 1.7 of [4]. □

PROOF OF PROPOSITION 1.1(i) AND (ii). Fix some real $s$ and $x$ and an arbitrarily small $\varepsilon > 0$. Find $\gamma > 0$ so small that $\delta_{[s-\gamma,s]}(\gamma) < \varepsilon$, a rational $u \in [s-\gamma, s]$, and a rational $y \in [x - 3\varepsilon, x - 2\varepsilon]$. It follows from Proposition 2.1 that $X_s^{u,y} \in (x - 4\varepsilon, x)$. Since $\varepsilon > 0$ is arbitrarily small, it follows from definition (1.5) that $X_s^{s,x-} = x$. Similarly, $X_s^{s,x+} = x$.

Lemma 2.2 and definitions (1.5) and (1.7) easily imply that a.s. for all $s, x \in \mathbb{R}$ and $t \geq s$ simultaneously, $X_t^{s,x+} \geq X_t^{s,x-}$.

It is an elementary fact that if all elements of an arbitrary family of functions have moduli of continuity bounded by $2\delta_{[a,b]}(\cdot)$, then so does the supremum of the functions in the family. This, definitions (1.5) and (1.7) and the remark following Proposition 2.1 imply that all processes $X_t^{s,x-}$ and $X_t^{s,x+}$ are Hölder continuous, for all real $s$ and $x$ simultaneously. □

LEMMA 2.3. (i) *For any rational times $s_1 < s_2$, with probability 1, the range of $\mathbb{Q} \ni x \to X_{s_2}^{s_1,x}$ consists of two semi-infinite "intervals" $(-\infty, y_1] \cap \mathbb{Q}$ and $[y_2, \infty) \cap \mathbb{Q}$, and a countable set $S \subset (y_1, y_2)$ which does not have accumulation points inside $(y_1, y_2)$. If $x \in \mathbb{Q}$ and $X_{s_2}^{s_1,x} \notin (y_1, y_2)$, then $L_{s_2}^{s_1,x} = 0$.*

(ii) *Fix any rational times $s_1 < s_2$ and let $\Lambda = \{X_{s_1}^{u,x} : u \in \mathbb{Q}, u \leq s_1, x \in \mathbb{Q}\}$. With probability 1, $\{X_{s_2}^{u,x} : u \in \mathbb{Q}, u \leq s_1, x \in \mathbb{Q}\}$ consists of two semi-infinite "intervals" $(-\infty, y_1] \cap \Lambda$ and $[y_2, \infty) \cap \Lambda$, and a countable set $S \subset (y_1, y_2)$ which does not have accumulation points inside $(y_1, y_2)$. If $u, x \in \mathbb{Q}$, $u \leq s_1$, and $X_{s_2}^{u,x} \notin (y_1, y_2)$, then $L_{s_2}^{u,x} - L_{s_1}^{u,x} = 0$.*

PROOF. (i) Standard arguments can be used to derive (i) from Theorem 1.2 of [4]. (ii) The argument proving (i) does not depend on the assumption that $x$'s are rational numbers but on the fact that rationals are countable. Since $\Lambda$ is countable, the same argument applies. □

LEMMA 2.4. *With probability 1, $s$ is not a bifurcation time for any $s$ which is a local extremum of $B_t$.*

PROOF. Consider the stochastic differential equation (1.3) with $s = 0$, rational $x \neq 0$, driven by a three-dimensional Bessel process $B_t$ in place of the Brownian motion. Recall that the path properties of the three-dimensional Bessel process are the same as those of the Brownian motion on any fixed time interval $[s_1, s_2]$, with $0 < s_1 < s_2 < \infty$, by the Cameron–Martin–Girsanov formula ([7], Section 3.5). For every $x \neq 0$, there exists (random) $s_3 > 0$ such that $B_t \neq -x$ for all $t \leq s_3$. Hence, we have strong existence and uniqueness for solutions to (1.3) with $s = 0$, simultaneously for all rational $x \neq 0$, driven by a three-dimensional Bessel process $B_t$, on every interval $[s_1, s_2] \subset (0, \infty)$ and, in fact, on the whole interval $[0, \infty)$.

Since the three-dimensional Bessel process is transient, it is easy to see that $P(L_\infty^{0,-1} < \infty) = p > 0$. Hence, $P(L_\infty^{0,-1} \leq b) = p_1 > 0$ for some $b < \infty$. By scaling, $P(L_\infty^{0,-x} \leq xb) = p_1$ for rational $x > 0$. Hence, for every $\varepsilon > 0$,

$$P\left(\bigcup_{x \in \mathbb{Q}, x \in (0,\varepsilon)} \{L_\infty^{0,-x} \leq xb\}\right) \geq p_1,$$

and so

$$P\left(\bigcap_{\varepsilon > 0} \bigcup_{x \in \mathbb{Q}, x \in (0,\varepsilon)} \{L_\infty^{0,-x} \leq xb\}\right) \geq p_1.$$

This and the definition of $X_t^{0,x-}$ and $L_t^{0,x-}$ show that $L_t^{0,0-} = 0$ for all $t \geq 0$ with probability greater than or equal to $p_1$. The event $A = \bigcup_{s>0} \{L_t^{0,0-} = 0, t \in [0, s]\}$ belongs to the germ $\sigma$-field $\mathcal{F}_{0+}$ and its probability is bounded below by $p_1 > 0$ so $P(A) = 1$, by Blumenthal's 0–1 law. If $A$ occurs, we must have $L_\infty^{0,0-} = 0$, because the three-dimensional Bessel process never returns to 0 a.s. This proves that $P(L_\infty^{0,0-} = 0) = 1$ a.s.

Now we go back to solutions of (1.3) driven by a Brownian motion $B_t$. Suppose $\beta \leq 0$, consider any rational numbers $0 \leq r_1 < r_2 < \infty$ and let $s$ denote the unique time when $B_t$ attains its minimum on $[r_1, r_2]$. Note that $s < r_2$ a.s. It is well known that $\{B_{t+s} - B_s, t \in [0, r_2 - s]\}$ has the same path properties as the three-dimensional Bessel process (this follows, e.g., from Williams' decomposition, see [15], Section VII.4). Hence, $L_t^{s,0-} = 0$ for $t \in [0, r_2 - s]$. We obviously have $L_t^{s,0+} = 0$ for $t \in [0, r_2 - s]$, so $s$ is not a bifurcation time. It is easy to see that $s$ is not a bifurcation time when $\beta > 0$. Every local minimum of $B_t$ is the global minimum over some interval $[r_1, r_2]$ with rational endpoints, so our argument holds for all local minima simultaneously. The local maxima can be dealt with in an analogous way. □

PROOF OF PROPOSITION 1.1(iii). Fix arbitrary real $s$ and $x$. Let $s_1$ be the smallest time greater than or equal to $s$ with the property that $B_{s_1} = -x + B_s$ and $B_{s_1}$ is not a local extremum of $B_t$. All local extrema of $B_t$ occur at different levels so there is at most one local extremum $s_2 \in [s, s_1)$ with the property that $B_{s_2} = -x + B_s$.

Consider $X_t^{s,x-}$ and first suppose that $L_{s_1}^{s,x-} = 0$. Then, clearly, $X_t^{s,x-} = x + B_t - B_s$ and $L_t^{s,x-}$ satisfy (1.3) on $[s, s_1]$. If there is no $s_2$ as described above, then obviously $L_t^{s,x-}$ satisfies (1.2) on $[s, s_1]$. Next assume that there exists a unique extremum $s_2 \in [s, s_1)$ with $B_{s_2} = -x + B_s$. By Trotter's theorem on the joint continuity of Brownian local time ([8], Section 5.1), the local time of $B_t$ at the level $-x + B_s$ does not increase between times $s$ and $s_1$. Hence, again $L_t^{s,x-}$ satisfies (1.2) on $[s, s_1]$.

Next we will keep the assumption that $L_{s_1}^{s,x-} = 0$ and consider $X_t^{s,x-}$ for $t > s_1$. Fix any rational $s_3 > s_1$. Recall that $s_1$ is not a local extremum. Since Brownian motion does not have points of increase, the range $\{B_t, t \in [s_1, s_3]\}$ contains an interval of length $\varepsilon > 0$ centered at $B_{s_1}$. It follows from Lemma 2.3 that

$$\{X_{s_3}^{u,y} : X_{s_1}^{u,y} \in [x + B_{s_1} - B_s - \varepsilon/2, x + B_{s_1} - B_s + \varepsilon/2], u, y \in \mathbb{Q}, u \leq s\}$$

is a finite set. This and the definition of $X_t^{s,x-}$ and $L_t^{s,x-}$ imply that $X_t^{s,x-} = X_t^{u,y}$ and $L_t^{s,x-} = L_t^{u,y}$ for some $u, y \in \mathbb{Q}$ and all $t \geq s_3$. Hence, $X_t^{s,x-}$ and $L_t^{s,x-}$ satisfy (1.2) and (1.3) on $[s_3, \infty)$. Since $s_3 > s_1$ is arbitrarily close to $s_1$, the conditions (1.2) and (1.3) are satisfied on $(s_1, \infty)$. We will now consider the case when $L_{s_1}^{s,x-} > 0$. It is easy to see that if there is no $s_2$ as defined at the beginning of the proof, then we must have $L_{s_1}^{s,x-} = 0$. Similarly, it is easy to see that $L_{s_2}^{s,x-} = 0$. Moreover, $\inf\{t : L_t^{s,x-} > 0\} = s_2$, because otherwise $L_t^{s,x-}$ would have been zero on the whole interval $[s, s_1]$.

It is easy to see that $s_2$ cannot be a local maximum of $B_t$. If it is, the Brownian motion $B_t$ has to stay below $-x + B_s$ for $t \in [s, s_2)$ in view of the definition of $s_2$. It follows that for any $t_1 \in (s, s_1)$, there exists $z_1 < x$ such that for $z \in (z_1, x)$, the Brownian motion $B_t$ does not hit $-z + B_s$ in the interval $[s, t_1)$. This implies that $L_{s_1}^{s,x-} = 0$, a contradiction.

We will now assume that $s_2$ is a local minimum of $B_t$. If $\beta > 0$, the definition of $X_t^{s,x-}$ implies that it is the sum of $x + B_t - B_s$ and a nondecreasing process. Hence, it stays above 0 on some interval $(s_2, t_2]$ with $t_2 > s_2$. This and the definition of $X_t^{s,x-}$ imply that for every $t_3 \in (s_2, t_2)$, there exists $z_2 < x$ such that for all $u, z \in \mathbb{Q}$ with $u < s$ and $X_s^{u,z} \in (z_2, x)$, the process $X_t^{u,z}$ stays above 0 on the interval $[t_3, t_2]$. Hence, $L_t^{u,z}$ does not increase on this interval and so the same can be said about $L_t^{s,x-}$. Since $t_3$ is arbitrarily close to $s_2$, we see that $L_t^{s,x-}$ does not increase on $[s_2, t_2]$—this contradicts the fact that $\inf\{t : L_t^{s,x-} > 0\} = s_2$.

Let us assume that $\beta < 0$. Suppose that for some $s_4 > s_2$ and all $t \in [s_2, s_4]$, we have $-\beta L_t^{s,x-} = \inf_{u \in [t,s_4]} B_u - B_{s_2}$. Let $s_5$ be the minimum of $B_t$ on some interval $(s_6, s_4)$ with the property that $s_5 \neq s_6, s_4$ and $s_6 \in (s_2, s_4)$. Then $-\beta L_t^{s_5, X_{s_5}^{s,x--}} \geq -\beta(L_t^{s,x-} - L_{s_5}^{s,x-})$ for $t \in [s_5, s_4]$ and $L_{s_4}^{s_5, X_{s_5}^{s,x-+}} = 0$. Thus, $s_5$ is a bifurcation point, but this contradicts Lemma 2.4. We conclude that in every right neighborhood of $s_2$ there exists $t$ with $-\beta L_t^{s,x-} > \inf_{u \in [t,s_4]} B_u - B_{s_2}$.

Consider $s_7 > s_2$, arbitrarily close to $s_2$, and find $s_8 \in (s_2, s_7)$ satisfying $-\beta L_{s_8}^{s,x-} > \inf_{u \in [s_8, s_7]} B_u - B_{s_2}$. Let $s_9 > s_8$ be the smallest time such that $B_{s_9} = B_{s_2} - \beta L_{s_8}^{s,x-}$ and $s_9$ is not a local extremum of $B_t$. Then we can repeat the argument applied above to $s_1$ and $s_3$ to see that $X_t^{s,x-} = X_t^{u,y}$ and $L_t^{s,x-} = L_t^{u,y}$ for some $u, y \in \mathbb{Q}$ and all $t \geq s_{10}$, for every $s_{10} > s_9$. By the uniform continuity of all processes $L$ and $X$ [Proposition 1.1(ii)], $s_9 \to s_2$ as $s_7 \to s_2$, so (1.2) and (1.3) hold for $X_t^{s,x-}$ and $L_t^{s,x-}$ on $[s_2, s_1]$.

The same proof applies to $X_t^{s,x+}$ and $L_t^{s,x+}$ by symmetry. □

LEMMA 2.5. *With probability 1, for all $s_1, s_2, x_1, x_2 \in \mathbb{R}$ simultaneously,*

$$X_t^{s_1, x_1-} = X_t^{s_1, x_1+} = X_t^{s_2, x_2-} = X_t^{s_2, x_2+}$$

*for $t \geq t_1$, where $t_1 < \infty$ depends on $s_1, s_2, x_1, x_2$.*

PROOF. Consider any rational $s_3 < \min(s_1, s_2)$ and note that, for sufficiently large $K < \infty$ and all rationals $x$ with $|x| > K$, $L_{s_1}^{s_3, x} = L_{s_2}^{s_3, x} = 0$. Choose rational $x_3$ and $x_4$ with the property that $X_{s_1}^{s_3, x_3} < x_1$, $X_{s_2}^{s_3, x_3} < x_2$, $X_{s_1}^{s_3, x_4} > x_1$, and $X_{s_2}^{s_3, x_4} > x_2$. By a theorem from [1], with probability 1, $X_t^{s_3, x_3} = X_t^{s_3, x_4}$ for $t \geq t_1$, where $t_1 = t_1(s_3, x_3, x_4) < \infty$. Note that this claim holds for all rational $s_3, x_3$ and $x_4$ simultaneously, a.s. The lemma now follows from $X_{t_1}^{s_3, x_3} = X_{t_1}^{s_3, x_4}$ and the definitions of $X_t^{s_1, x_1-}, X_t^{s_1, x_1+}, X_t^{s_2, x_2-}$ and $X_t^{s_2, x_2+}$. □

LEMMA 2.6. *With probability 1, for all pairs $s, s_1 \in \mathbb{R}$ simultaneously, if $s < s_1$, then the set $\{x \in \mathbb{R} : X_{s_1}^{s,x-} \neq X_{s_1}^{s,x+}\}$ is countable.*

PROOF. This follows easily from Lemma 2.3. □

LEMMA 2.7. *With probability 1, for any $s, x \in \mathbb{R}$ and $t_1 > s$, there exist $s_1, x_1 \in \mathbb{Q}$, such that $X_t^{s,x-} = X_t^{s_1, x_1}$ for $t \geq t_1$, unless $L_{t_1}^{s,x} = 0$, and a similar statement holds for $X_t^{s,x+}$.*

PROOF. This follows from the argument given in the proof of Proposition 1.1(iii). □

Suppose $s_1, s_2 \in \mathbb{Q}$, $s_1 < s_2$, and consider points $y_1$ and $y_2$ in $\{X^{s_1,x}_{s_2}, x \in \mathbb{Q}\}$, such that $y_1 < y_2$ and $(y_1, y_2) \cap \{X^{s_1,x}_{s_2}, x \in \mathbb{Q}\} = \varnothing$. Points $y_1$ and $y_2$ with these properties exist by the results of [4] (see Lemma 2.3 above). The results of [4] show in addition that the sets $\{x \in \mathbb{Q} : X^{s_1,x}_{s_2} = y_1\}$ and $\{x \in \mathbb{Q} : X^{s_1,x}_{s_2} = y_2\}$ are intervals in $\mathbb{Q}$ with a common endpoint $z$. It follows easily that $y_1 = X^{s_1,z-}_{s_2} < X^{s_1,z+}_{s_2} = y_2$. Hence, we see that we cannot have strong uniqueness of solutions to (1.3) simultaneously for all $s, x \in \mathbb{R}$, and so bifurcation times exist. Note that for a fixed $s_1$, typically there are many $z$'s with $X^{s_1,z-}_{s_2} < X^{s_1,z+}_{s_2}$. For all $z \neq 0$, in this family we have $X^{s_1,z-}_t = X^{s_1,z+}_t$ for $t \in [s_1, s_3]$ and some $s_3 = s_3(z) > s_1$.

Recall that if $\{(X^{s,x-}_t, X^{s,x+}_t), t \in [s,u]\}$ is a lens, then $x = X^{s,x-}_u = X^{s,x+}_u = 0$.

PROOF OF THEOREM 1.3(i). Consider any lens $\{(X^{s,0-}_t, X^{s,0+}_t), t \in [s,u]\}$ and note that by Lemma 2.7, for some $s_1, s_2, x_1, x_2 \in \mathbb{Q}$ and $u_1 \in (s,u)$, we have $X^{s,0-}_t = X^{s_1,x_1}_t$ and $X^{s,0+}_t = X^{s_2,x_2}_t$ for all $t \in [u_1, u]$. Hence, $(u, X^{s_1,x_1}_u) = (u, 0)$ is a point in space-time where the processes $X^{s_1,x_1}_\cdot$ and $X^{s_2,x_2}_\cdot$ coalesce, for some rational $s_1, s_2, x_1, x_2 \in \mathbb{Q}$. The set of such points is countable. It is easy to see that it is also infinite. □

Note that in the next two lemmas $T_1$ denotes a bifurcation time.

LEMMA 2.8. *The following holds with probability 1 for all $s, t_1, x \in \mathbb{R}$ simultaneously. Suppose that $s < t_1$ and $X^{s,x-}_{t_1} < X^{s,x+}_{t_1}$. Let $T_1 = \inf\{t \geq s : X^{s,x-}_t \neq X^{s,x+}_t\}$ and $T_2 = \sup\{t \geq s : L^{s,x-}_t \vee L^{s,x+}_t = 0\}$. Then $T_1 = T_2$.*

PROOF. Clearly, $T_1 \geq T_2$. Suppose that $T_1 > T_2$, note that $t_1 \geq T_1$ and let $t_2 \in (T_2, T_1)$. Then, by Lemma 2.7, there exist $s_1, s_2, x_1, x_2 \in \mathbb{Q}$ such that $X^{s,x-}_t = X^{s_1,x_1}_t$ and $X^{s,x+}_t = X^{s_2,x_2}_t$ for $t \geq t_2$. By the strong uniqueness of solutions to (1.3) for rational $s$ and $x$, $X^{s,x-}_t = X^{s_1,x_1}_t = X^{s_2,x_2}_t = X^{s,x+}_t$ for $t \geq t_2$. This contradicts the assumption that $X^{s,x-}_{t_1} < X^{s,x+}_{t_1}$ and so $T_1 = T_2$. □

LEMMA 2.9. *With probability 1, the following holds for all $s, x \in \mathbb{R}$. Suppose that $X^{s,x-}_{t_1} < X^{s,x+}_{t_1}$ for some $t_1 > s$ and $x \in \mathbb{R}$. Let $T_1 = \inf\{t \geq s : X^{s,x-}_t < X^{s,x+}_t\}$. If $x \neq 0$, then $T_1 > s$ and $T_1$ is not a semi-flat bifurcation time.*

PROOF. It is obvious that $x \neq 0$ implies $T_1 > s$. By Lemma 2.4, the bifurcation time $T_1$ is not a local extremum of $B_t$. Let $y = B_s - x$. Note that $B_{T_1} = y$ and either $B_t > y$ for all $t \in [s, T_1)$ or $B_t < y$ for $t \in [s, T_1)$. Since Brownian motion does not have points of increase (or decrease) and $T_1$ is not a local extremum, it follows that $B_t$ crosses the level $y$ infinitely often in every interval $(T_1, T_1 + \delta)$, $\delta > 0$. By Lemma 2.8, $L_{T_1}^{s,x-} = L_{T_1}^{s,x+} = 0$. Suppose that $T_1$ is a semi-flat bifurcation time. Then for some $s_1 > T_1$, either $L_{s_1}^{s,x-} = 0$ or $L_{s_1}^{s,x+} = 0$. Assume without loss of generality that $L_{s_1}^{s,x-} = 0$. Since the Brownian motion $B_t$ crosses the level $y$ repeatedly between times $T_1$ and $s_1$, it accumulates some local time at this level (by Trotter and Ray–Knight theorems, [8], Sections 5.1 and 5.3), and so the process $X_{s_1}^{s,x-}$ accumulates some local time at 0. This contradicts the assumption that $L_{s_1}^{s,x-} = 0$. □

PROOF OF THEOREM 1.3(ii). Recall that we call a bifurcation time *ordinary* if it is not semi-flat. First, we are going to show that if $s$ is an ordinary bifurcation time and $s_1 > s$ is such that $X_{s_1}^{s,0-} \neq X_{s_1}^{s,0+}$, then for every $\varepsilon > 0$, there exists an ordinary bifurcation time $s_2 \in (s, s_1 \wedge (s + \varepsilon/2))$ such that $X_{s_2}^{s,0-} < 0 < X_{s_2}^{s,0+}$, $X_{s_1}^{s_2,0-} \neq X_{s_1}^{s_2,0+}$, $L_{s+2(s_2-s)}^{s,0-} < L_{s+\varepsilon}^{s,0-}$, and $L_{s+2(s_2-s)}^{s,0+} < L_{s+\varepsilon}^{s,0+}$. Since $s$ is an ordinary bifurcation time, we can find $s_3 \in (s, s_1 \wedge (s + \varepsilon/2))$ such that both processes $X_{\cdot}^{s,0-}$ and $X_{\cdot}^{s,0+}$ cross 0 in the interval $(s_3, s_1 \wedge (s + \varepsilon/2))$. This and Lemma 2.3 imply that the set $\{X_{s_1}^{t,y}, t \leq s_3, X_{s_3}^{t,y} \in [X_{s_3}^{s,0-}, X_{s_3}^{s,0+}], t, y \in \mathbb{Q}\}$ is finite; let $x_1$ be the second largest element of this set and denote $x_2 = X_{s_1}^{s,0+}$, that is, the largest element of the set. For $t < s_1$, let $\Gamma_t = \inf\{y \in \mathbb{Q} : X_{s_1}^{t,y} = X_{s_1}^{s,0+}\}$ and note that $\Gamma_t$ is not constantly equal to 0 on any interval $(s, s + \delta)$ because $s$ is an ordinary bifurcation time. Let $s_4 \in (s, s + (s_3 - s)/2)$ be so close to $s$ that $X_{\cdot}^{s,x-}$ and $X_{\cdot}^{s,x+}$ cross 0 in the interval $(s_4, s + (s_3 - s)/2)$, and $\Gamma_{s_4} \neq 0$. Note that $X_{s_1}^{s_4, \Gamma_{s_4}-} \leq x_1 < x_2 \leq X_{s_1}^{s_4, \Gamma_{s_4}+}$ and let $s_2 = \inf\{t \geq s_4 : X_t^{s_4, \Gamma_{s_4}-} < X_t^{s_4, \Gamma_{s_4}+}\}$. By Lemma 2.9, $s_2$ is an ordinary bifurcation time. Note that $s_2 \leq s + (s_3 - s)/2 \leq s_1 \wedge (s + \varepsilon/2)$ because $X_{\cdot}^{s,0-}$ and $X_{\cdot}^{s,0+}$ cross 0 in the interval $(s_4, s + (s_3 - s)/2)$. Moreover, $X_{s_1}^{s_2, 0-} = X_{s_1}^{s_4, \Gamma_{s_4}-} \leq x_1 < x_2 \leq X_{s_1}^{s_4, \Gamma_{s_4}+} = X_{s_1}^{s_2, 0+}$. We have $L_{s+2(s_2-s)}^{s,0-} < L_{s+\varepsilon}^{s,0-}$, and $L_{s+2(s_2-s)}^{s,0+} < L_{s+\varepsilon}^{s,0+}$ because $X_{\cdot}^{s,0-}$ and $X_{\cdot}^{s,0+}$ cross 0 in the interval $(s_3, s_1 \wedge (s + \varepsilon/2))$. It is easy to see that $X_{s_2}^{s,0-} \leq 0 \leq X_{s_2}^{s,0+}$. The inequalities are, in fact, sharp because on every interval $(s_5, \infty)$ with $s_5 > s$, $X_{\cdot}^{s,0-}$ agrees with some $X_{\cdot}^{u,z}$ with rational $u$ and $z$, and the last process does not pass through bifurcation points, a.s., and the same holds for $X_{\cdot}^{s,0+}$. This completes the proof of our claim.

Recall that $s$ is an ordinary bifurcation time and $s_1 > s$ is such that $X_{s_1}^{s,0-} \neq X_{s_1}^{s,0+}$. We will construct a family of ordinary bifurcation times $u$ with $u \in [s, ((s+s_1)/2) \wedge s_3)$ in the following inductive way. Start with an ordinary bifurcation time $s_0$ such that $s_0 \in (s, (s+(s+s_1)/2) \wedge s_3)$, $X_{s_0}^{s,0-} < 0 < X_{s_0}^{s,0+}$, $X_{s_1}^{s_0, 0-} \leq x_1$, $X_{s_1}^{s_0, 0+} \geq x_2$, $L_{s+2(s_0-s)}^{s,0-} < L_{s+1}^{s,0-}$, and $L_{s+2(s_0-s)}^{s,0+} < L_{s+1}^{s,0+}$. Then find ordinary bifurcation times $s_{00}$ and $s_{01}$ such that $s_0 <$

$s_{00} < s_{01} < (s+2^{-1}) \wedge s_3$, $X^{s_0,0-}_{s_{00}} < 0 < X^{s_0,0+}_{s_{00}}$, $X^{s_{00},0-}_{s_1} \leq x_1$, $X^{s_{00},0+}_{s_1} \geq x_2$, $L^{s_0,0-}_{s_0+2(s_{00}-s_0)} < L^{s_0,0-}_{s+2^{-1}}$, $L^{s_0,0+}_{s_0+2(s_{00}-s_0)} < L^{s_0,0+}_{s+2^{-1}}$, $X^{s_0,0-}_{s_{01}} < 0 < X^{s_0,0+}_{s_{01}}$, $X^{s_{01},0-}_{s_1} \leq x_1$, $X^{s_{01},0+}_{s_1} \geq x_2$, $L^{s_0,0-}_{s_0+2(s_{01}-s_0)} < L^{s_0,0-}_{s+2^{-1}}$, and $L^{s_0,0+}_{s_0+2(s_{01}-s_0)} < L^{s_0,0+}_{s+2^{-1}}$.

We can find inductively ordinary bifurcation times $s_{0k_1k_2...k_n}$ for all $n \geq 2$, with $k_j = 0, 1$, with the following properties. Suppose that $s_{0k_1k_2...k_n}$ have been defined for $n \leq m$ and let $\delta_m$ be the minimum of distances between distinct elements of $\{s_{0k_1k_2...k_n}, n \leq m, k_j = 0, 1\}$. Then for any $k_j = 0, 1$, $j = 1, 2, \ldots, m$, find ordinary bifurcation times $s_{0k_1k_2...k_m0}$ and $s_{0k_1k_2...k_m1}$ with

$$s_{0k_1k_2...k_m} < s_{0k_1k_2...k_m0} < s_{0k_1k_2...k_m1}$$
$$< (s_{0k_1k_2...k_m} + (\delta_m/10) \wedge 2^{-m+1}) \wedge s_3,$$

$$X^{s_{0k_1k_2...k_m},0-}_{s_{0k_1k_2...k_m0}} < 0 < X^{s_{0k_1k_2...k_m},0+}_{s_{0k_1k_2...k_m0}},$$

$$X^{s_{0k_1k_2...k_m0},0-}_{s_1} \leq x_1,$$

$$X^{s_{0k_1k_2...k_m0},0+}_{s_1} \geq x_2,$$

$$v_{m0} = s_{0k_1k_2...k_m} + 2(s_{0k_1k_2...k_m0} - s_{0k_1k_2...k_m}),$$

$$L^{s_{0k_1k_2...k_m},0-}_{v_{m0}} < L^{s_{0k_1k_2...k_m},0-}_{s+2^{-m}},$$

$$L^{s_{0k_1k_2...k_m},0+}_{v_{m0}} < L^{s_{0k_1k_2...k_m},0+}_{s+2^{-m}},$$

(2.1) $$X^{s_{0k_1k_2...k_m},0-}_{s_{0k_1k_2...k_m1}} < 0 < X^{s_{0k_1k_2...k_m},0+}_{s_{0k_1k_2...k_m1}},$$

$$X^{s_{0k_1k_2...k_m1},0-}_{s_1} \leq x_1,$$

$$X^{s_{0k_1k_2...k_m1},0+}_{s_1} \geq x_2,$$

$$v_{m1} = s_{0k_1k_2...k_m} + 2(s_{0k_1k_2...k_m1} - s_{0k_1k_2...k_m}),$$

$$L^{s_{0k_1k_2...k_m},0-}_{v_{m1}} < L^{s_{0k_1k_2...k_m},0-}_{s+2^{-m}},$$

$$L^{s_{0k_1k_2...k_m},0+}_{v_{m1}} < L^{s_{0k_1k_2...k_m},0+}_{s+2^{-m}}.$$

It follows from (2.1) that the sequence $(s_{0k_1}, s_{0k_1k_2}, s_{0k_1k_2k_3}, \ldots)$ converges for any choice of $0, k_1, k_2, k_3, \ldots$. The family $\mathcal{S}$ of limit points is uncountable. We will show that every element of $\mathcal{S}$ is a bifurcation time. Fix any $u \in \mathcal{S}$ and any $u_1 \in (u, s_1)$. Find $0, k_1, k_2, k_3, \ldots$ such that $s_{0,k_1,k_2,k_3,\ldots,k_m} \uparrow u$. By construction, $X^{s_{0k_1k_2...k_m},0-}_{\cdot} \leq X^{s_{0k_1k_2...k_mk_{m+1}},0-}_{\cdot}$ and $X^{s_{0k_1k_2...k_m},0+}_{\cdot} \geq X^{s_{0k_1k_2...k_mk_{m+1}},0+}_{\cdot}$ on the interval $(s_{0k_1k_2...k_mk_{m+1}}, \infty)$. Passing to the limit and using the fact that $X^{s_{0k_1k_2...k_m},0-}_{s_{0k_1k_2...k_m0}} < 0 < X^{s_{0k_1k_2...k_m},0+}_{s_{0k_1k_2...k_m0}}$ and $X^{s_{0k_1k_2...k_m},0-}_{s_{0k_1k_2...k_m1}} <$ $0 < X^{s_0 k_1 k_2 \ldots k_m, 0+}_{s_0 k_1 k_2 \ldots k_m 1}$, we see that $X^{s_0 k_1 k_2 \ldots k_m, 0-}_u \leq 0 \leq X^{s_0 k_1 k_2 \ldots k_m, 0+}_u$. For sufficiently large $m$, we have $L^{s_0 k_1 k_2 \ldots k_m, 0-}_u < L^{s_0 k_1 k_2 \ldots k_m, 0-}_{u_1}$ and $L^{s_0 k_1 k_2 \ldots k_m, 0+}_u < L^{s_0 k_1 k_2 \ldots k_m, 0+}_{u_1}$. This implies that $L^{u, 0+}_\cdot$ and $L^{u, 0-}_\cdot$ increase on the interval $(u, u_1)$. Hence, for every $u_2 > u$, $X^{u, 0+}_\cdot$ agrees with some $X^{v, x}_\cdot$ with rational $v$ and $x$ on the interval $(u_2, \infty)$, and the same remark applies to $X^{u, 0-}_\cdot$. Since $X^{v, x}_\cdot$ does not pass through any bifurcation points, either $u$ is a bifurcation time or $X^{u, 0-}_t = X^{u, 0+}_t$ for all $t \geq u$. Our construction implies that $X^{u, 0-}_{s_1} \leq x_1 < x_2 \leq X^{u, 0+}_{s_1}$ so $u$ must be a bifurcation time. $\square$

We will use excursion theory in the next proof and the proof of Theorem 1.6. Various accounts of excursion theory may be found in [2, 15], Chapter XII, and [16], Chapter 8. Some of the most relevant material is contained in [13]. We will use Proposition 4.1 and Theorem 5.1 of [3].

PROOF OF THEOREM 1.4. (i) We will use the method of Watanabe [19]. First we will construct a family of "excursions" whose starting points are semi-flat bifurcation points. Then we will assemble these excursion into a Brownian motion, as in [19].

Assume without loss of generality that $\beta \in (-1, -\frac{1}{3})$ and fix a large $K < \infty$ whose value will be specified later. Let $\widetilde{B}^1_t, \widetilde{B}^2_t, \ldots$ be independent Brownian motions starting from 0 and let $B^k_t = \widetilde{B}^k_t - \beta 2^{-k}$. Define processes $X^k_t$ by equations analogous to (1.3) and (1.2):

$$X^k_t = B^k_t + \beta L^k_t, \qquad t \geq 0,$$

$$L^k_t = 2^{-k} + \lim_{\varepsilon \to 0} \frac{1}{2\varepsilon} \int_0^t \mathbf{1}_{(-\varepsilon, \varepsilon)}(X^k_s) \, ds.$$

Let

$$T_1 = \inf\{t > 0 : B^1_t = 0 \text{ or } X^1_t \geq -K\beta L^1_t\},$$

$$T_k = \inf\{t > 0 : B^k_t = 0 \text{ or } L^k_t = 2^{-k+1} \text{ or } X^k_t \geq -K\beta L^k_t\}, \qquad k \geq 2,$$

$$S_{k,n} = T_k + T_{k-1} + \cdots + T_{n+1}, \qquad n = 0, 1, \ldots, k-1,$$

$$Y^k_t = \begin{cases} X^k_t, & \text{for } t \in [0, T_k); \\ X^n_{t - S_{k,n}}, & \text{for } t \in [S_{k,n}, S_{k,n-1}), \, 1 \leq n \leq k-1, \\ & \text{if } B^m_{T_m} > 0 \text{ and } X^m_{T_m} < -K\beta L^m_{T_m} \text{ for all } m < n; \\ 0, & \text{for } t \in [S_{k,n}, S_{k,n-1}), \, 1 \leq n \leq k-1, \\ & \text{if } B^m_{T_m} = 0 \text{ or } X^m_{T_m} = -K\beta L^m_{T_m} \text{ for some } m < n; \\ 0, & \text{for } t \geq S_{k,0}, \end{cases}$$

$$A_t^k = \begin{cases} B_t^k, & \text{for } t \in [0, T_k); \\ B_{t-S_{k,n}}^n, & \text{for } t \in [S_{k,n}, S_{k,n-1}), \ 1 \leq n \leq k-1, \\ & \text{if } B_{T_m}^m > 0 \text{ and } X_{T_m}^m < -K\beta L_{T_m}^m \text{ for all } m < n; \\ 0, & \text{for } t \in [S_{k,n}, S_{k,n-1}), \ 1 \leq n \leq k-1, \\ & \text{if } B_{T_m}^m = 0 \text{ or } X_{T_m}^m = -K\beta L_{T_m}^m \text{ for some } m < n; \\ 0, & \text{for } t \geq S_{k,0}. \end{cases}$$

In other words, $Y_t^k$ is a process assembled from $X_t^k, X_t^{k-1}, \ldots, X_t^1$ and $A_t^k$ is a Brownian motion assembled from $B_t^k, B_t^{k-1}, \ldots, B_t^1$. The processes $Y_t^k$ and $A_t^k$ are sent to 0 at the time

$$S_k = \inf\{t > 0 : A_t^k = 0 \text{ or } Y_t^k \geq -K\beta \widehat{L}_t^k\},$$

where

$$\widehat{L}_t^k = 2^{-k} + \lim_{\varepsilon \to 0} \frac{1}{2\varepsilon} \int_0^t \mathbf{1}_{(-\varepsilon, \varepsilon)}(Y_s^k) \, ds.$$

It is elementary to check that the processes $Y_t^k$ and $A_t^k$ satisfy the equation analogous to that for $X_t^k$ and $B_t^k$:

$$Y_t^k = A_t^k + \beta \widehat{L}_t^k, \qquad 0 \leq t \leq S_k.$$

Let $Q^k$ denote the distribution of $A_t^k$ and let $\alpha = \frac{\beta-1}{2\beta} - \frac{1+\beta}{2K\beta}$. We will argue that $\lim_{k \to \infty} 2^{\alpha k} Q^k$ exists and defines an excursion law $H$ for Brownian motion.

Let $U_k = \inf\{t > 0 : \widehat{L}_t^k \geq 2^{-k+1}\}$. The event $\{U_k < S_k\}$ is the same as the first excursion of $Y_t^k$ above 0 of height greater than $-K\beta \widehat{L}_t^k$ and the first excursion of $Y_t^k$ below 0 which reaches the level $\beta \widehat{L}_t^k$ occuring after $\widehat{L}_t^k$ increases by $2^{-k}$. According to the excursion theory, on the time scale corresponding to the local time $\widehat{L}_t^k$, the point process of arrivals of excursions of $Y_t^k$ above 0 with height greater than $-K\beta \widehat{L}_t^k$ is a Poisson point process with the variable intensity, equal to $f_1(s) = -(1+\beta)/(2K\beta(s+2^{-k}))$. The intensity for the analogous point process of arrivals of excursions of $Y_t^k$ below 0 which hit $\beta \widehat{L}_t^k$ is equal to $f_2(s) = -(1-\beta)/(2\beta(s+2^{-k}))$. Hence, the probability that none of these excursions occurs before $\widehat{L}_t^k$ increases by $2^{-k}$ is equal to

$$\exp\left(-\int_0^{2^{-k}} \left(\frac{-(1+\beta)}{2K\beta(s+2^{-k})} + \frac{\beta-1}{2\beta(s+2^{-k})}\right) ds\right).$$

Elementary calculations show that this is equal to $2^{-\alpha}$. Let $U_k^n = \inf\{t > 0 : \widehat{L}_t^k \geq 2^{-n}\}$ for $n < k$. By induction and the strong Markov property,

(2.2) $$P(U_k^n < S_k) = 2^{-(k-n)\alpha}.$$

It follows that for every fixed $n$ and all $k > n$, the measures $2^{\alpha k} Q^k$ give the same mass to paths of $Y_t^k$ in the set $\{U_k^n < S_k\}$. It is clear from the construction of processes $Y_t^k$ that for a fixed $n$, the conditional distribution $F_k^n$ of $U_k^n$ given $\{U_k^n < S_k\}$ is nondecreasing in $k$, that is, the distribution $F_{k+1}^n$ is stochastically larger than $F_k^n$. If we show that the expectations of $F_k^n$ are uniformly bounded in $k$, that will prove that the distributions $F_k^n$ converge as $k \to \infty$.

Let $c_1$ be the expected lifetime of a Brownian excursion above 0 conditioned on not hitting level 1. It is well known that $c_1 < \infty$. By scaling, the expectation of excursion lifetime conditioned on not hitting level $a$ is equal to $c_1 a^2$. For $k > n$, the expectation of $F_k^n$ is equal to

$$(2.3) \qquad \int_{2^{-k}}^{2^{-n}} \left( \frac{1+\beta}{2} c_1 (K\beta s)^2 + \frac{1-\beta}{2} c_1 (\beta s)^2 \right) ds \leq c_2 2^{-3n},$$

so $F_k^n$'s converge as $k \to \infty$. This and the strong Markov property applied at $U_k^n$ imply that the distributions of $\{Y_t^k, t \geq U_k^n\}$ under $2^{\alpha k} Q^k$ converge as $k \to \infty$. Since $\{A_t^k, 0 \leq t \leq s\}$ is a function of $\{Y_t^k, 0 \leq t \leq s\}$, we have similar convergence for distributions of $A_t^k$'s. The integer $n$ is arbitrary, so we conclude that $H = \lim_{k \to \infty} 2^{\alpha k} Q^k$ exists.

It is clear from the definition that $H$ is a $\sigma$-finite measure which is the "distribution" of a process $\{A_t, 0 \leq t \leq S\}$ satisfying

$$(2.4) \qquad \begin{aligned} Y_t &= A_t + \beta \widehat{L}_t, \qquad 0 \leq t \leq S, \\ \widehat{L}_t &= \lim_{\varepsilon \to 0} \frac{1}{2\varepsilon} \int_0^t \mathbf{1}_{(-\varepsilon, \varepsilon)}(Y_s)\, ds, \\ S &= \inf\{t > 0 : A_t = 0 \text{ or } Y_t \geq -K\beta \widehat{L}_t\}, \end{aligned}$$

and such that for every $s > 0$, the distribution of $\{A_t, t \in [s, S]\}$ given $\{S > s\}$ is a Brownian motion stopped at $S$.

By (2.3), the $H$-expectation of $S \wedge \inf\{t > 0 : \widehat{L}_t = 1\}$ is bounded by $\sum_{n \geq 0} c_2 \times 2^{-3n} 2^{n\alpha} < \infty$, assuming $\alpha < 3$. Since the $H$-measure of $\{\inf\{t > 0 : \widehat{L}_t = 1\} < S\}$ is finite, the $H$-measure of $\{S > s\}$ is finite for every $s > 0$. It follows from the definition of $H$ and (2.2) that the $H$-measure of $\widehat{L}$-paths hitting level $2^{-n}$ is equal to $2^{n\alpha}$. This and the easy fact that $H$ has the same space-time scaling properties as the Brownian motion imply that $H(S > s) = c_3 s^{-\alpha/2}$. Hence, the $H$-expectation of $S$ on the set $\{S < 1\}$ is finite if $\alpha < 2$. Recall that $\beta \in (-1, -\frac{1}{3})$ and note that we can choose $K$ so large that $\alpha < 2$.

Now generate a Poisson point process of excursions on the product of $[0, \infty)$ and the space of stopped continuous paths with intensity given by the product of the Lebesgue measure and $H$. The excursions can be assembled into a Brownian path, as in [19], because the $H$-expectation of $S$ on the set $\{S < 1\}$ is finite. The starting points of constituent excursions in this Brownian path are semi-flat bifurcation points—this follows from (2.4) and the fact that $S > 0$ for $H$-almost every path.

(ii) This part of the proof is based on a classical covering argument.

Suppose without loss of generality that $\beta \in (-1/3, 0)$ and let $\gamma = (\beta - 1)/(2\beta)$. Then $\gamma > 2$. Let $S = \inf\{t > 0 : L_t^{0,0} = 1\}$ and $T_\varepsilon = \inf\{t > 0 : B_t = -\varepsilon\}$, for $\varepsilon > 0$. The event $\{S < T_\varepsilon\}$ is the same as that the first excursion of $X^{0,0}$ below 0, starting at a time $s$, which hits the level $-\varepsilon + \beta L_s^{0,0}$ occurring after time $S$. We calculate the probability of this event using excursion theory, as in part (i) of the proof,

$$(2.5) \quad P(S < T) = \exp\left(-\int_0^1 \frac{1-\beta}{2(\varepsilon - \beta s)}\,ds\right) = (1 - \beta/\varepsilon)^{(1-\beta)/2\beta} \leq c_1 \varepsilon^\gamma.$$

Let $S_{s,\varepsilon} = \inf\{t > 0 : L_t^{s,-\varepsilon} = 1\}$ and $T_{s,\varepsilon} = \inf\{t > 0 : X_t^{s,\varepsilon} = 0\}$, for $\varepsilon > 0$. Note that $T_{s,\varepsilon} = \inf\{t > 0 : B_t - B_s = -\varepsilon\}$. By the strong Markov property applied at $\inf\{t > 0 : X_t^{s,-\varepsilon} = 0\}$, (2.5), scaling and shift invariance of Brownian motion, we have

$$(2.6) \quad P(S_{s,\varepsilon} < T_{s,\varepsilon}) \leq c_2 \varepsilon^\gamma.$$

Fix some $n > 1$ and $\alpha \in (2, \gamma)$, let $\varepsilon = 1/n$ and $s_k = k/n$, for $k = 0, 1, \ldots, n$. Let $A_k$ be the event that there exists a semi-flat bifurcation time $u \in [s_{k-1}, s_k]$ with the following properties: $X_{s_k}^{u,0+}, X_{s_k}^{u,0-} \in [-\varepsilon^{1/\alpha}, \varepsilon^{1/\alpha}]$ and $\inf\{t > 0 : L_t^{u,0-} = 1\} < \inf\{t > 0 : X_t^{u,0+} = 0\}$. Then Lemma 2.2 easily implies that on $A_k$, $S_{s_k, \varepsilon^{1/\alpha}} < T_{s_k, \varepsilon^{1/\alpha}}$. By (2.6), $P(A_k) \leq c_2 \varepsilon^{\gamma/\alpha}$ and

$$(2.7) \quad P\left(\bigcup_{1 \leq k \leq n} A_k\right) \leq c_2 \varepsilon^{\gamma/\alpha - 1}.$$

The following standard estimate for Brownian motion,

$$(2.8) \quad P\left(\sup_{t \in [s_{k-1}, s_k]} |B_t - B_{s_{k-1}}| \geq \varepsilon^{1/\alpha}\right) \leq c_3 \varepsilon^{1/\alpha - 1/2} \exp(-\tfrac{1}{2}\varepsilon^{2/\alpha - 1}),$$

applies also to all skew Brownian motions driven by $B_t$ because of Proposition 2.1. Hence, if $C_k$ denotes the event that there exists a semi-flat bifurcation time $u \in [s_{k-1}, s_k]$ with $\inf\{t > 0 : L_t^{u,0-} = 1\} < \inf\{t > 0 : X_t^{u,0+} = 0\}$, then (2.7) and (2.8) yield

$$P\left(\bigcup_{1 \leq k \leq n} C_k\right) \leq c_2 \varepsilon^{\gamma/\alpha - 1} + \varepsilon^{-1} c_3 \varepsilon^{1/\alpha - 1/2} \exp(-\tfrac{1}{2}\varepsilon^{2/\alpha - 1}).$$

This goes to 0 as $\varepsilon \to 0$ so there are no semi-flat bifurcation times $u \in [0, 1]$ with $\inf\{t > 0 : L_t^{u,0-} = 1\} < \inf\{t > 0 : X_t^{u,0+} = 0\}$. An analogous argument shows that for any integer $m > 1$, there are no semi-flat bifurcation times $u \in [0, m]$ with $\inf\{t > 0 : L_t^{u,0-} = 1/m\} < \inf\{t > 0 : X_t^{u,0+} = 0\}$. Hence, there are no semi-flat bifurcation times if $|\beta| < \frac{1}{3}$. □

PROOF OF THEOREM 1.5(i). Note that the distribution of $L_t^{0,0}$ does not depend on $\beta$. Since for $\beta = 0$ this is the usual Brownian local time, $L_1^{0,0}$ has a continuous density, for any $\beta$. An easy argument based on the strong Markov property and scaling shows that for any $x$, the random variable $x + \beta L_1^{0,x}$ has a continuous density on $(x, \infty)$. This clearly implies that with probability 1, for all rational $x$ simultaneously, $x + \beta L_1^{0,x} \neq 0$. Hence, by Theorem 1.2 of [4], $\inf_{x \in \mathbb{Q}} |x + \beta L_1^{0,x}| > 0$, a.s. Let $\rho_t = \inf_{x \in \mathbb{Q}} |x + \beta L_t^{0,x}|$. By scaling, for $t > 0$, the distribution of $\rho_t$ is the same as that of $\sqrt{t}\rho_1$, and it is also the same as the distribution of $\inf_{x \in \mathbb{Q}} |x + \beta L_1^{1-t,x}|$. Hence, for every $a, p > 0$, $P(\inf_{x \in \mathbb{Q}} |x + \beta L_1^{1-t,x}| < a) < p$, for sufficiently large $t$. Note that the random set $\{x + \beta L_1^{t,x} : x \in \mathbb{Q}\}$ is increasing in $t$. It follows that for every $a, p > 0$, and sufficiently large $t$,

$$(2.9) \qquad P\left(\inf_{x \in \mathbb{Q}} \inf_{u > t} |x + \beta L_1^{1-u,x}| < a\right) < p.$$

Let $\Gamma = \{X_1^{0,x} : x \in \mathbb{Q}\}$, $y_1 = X_1^{0,0}$, $y_2 = \inf\{y \in \Gamma : y > y_1\}$ and $y_3 = \sup\{y \in \Gamma : y < y_1\}$. It follows from Theorem 1.2 of [4] that $y_1$ is an isolated point in $\Gamma$, so $y_3 < y_1 < y_2$, a.s. Let $\Lambda_t = \{x \in \mathbb{Q} : X_1^{t,x} = y_1\}$ and $s_0 = \inf\{t : \Lambda_t \neq \varnothing\}$. We will prove that $s_0 > -\infty$ and then we will show that three distinct solutions to (1.3) start at $(s_0, 0)$.

Suppose that $s_0 = -\infty$ with positive probability. Find $a, p_1 > 0$ such that $P(|y_1| < a/2, s_0 = -\infty) > p_1$. According to (2.9), we can find $t > 0$ so large that

$$P\left(\inf_{x \in \mathbb{Q}} \inf_{u > t} |x + \beta L_1^{1-u,x}| < a\right) < p_1/2.$$

Find $t_1 > t$ such that with probability greater than $1 - p_1/4$, there exists $t_2 \in (t, t_1)$, such that $B_{1-t_2} = B_1$. Then $X_1^{1-t_2, x} = x + \beta L_1^{1-t_2,x}$ for all rational $x$, and, hence, $s_0 \geq t_2 > t_1$. We see that $P(s_0 < t_1) < p_1/2 + p_1/4$, which contradicts the assumption that $P(|y_1| < a/2, s_0 = -\infty) > p_1$. It follows that $P(s_0 > -\infty) = 1$.

Let $\lambda_t^- = \inf \Lambda_t$, $\lambda_t^+ = \sup \Lambda_t$, $\gamma^- = \liminf_{t \downarrow s_0} \lambda_t^-$, and $\gamma^+ = \limsup_{t \downarrow s_0} \lambda_t^+$. We will first prove that $\gamma^- = \gamma^+$. Suppose that $\gamma^+ - \gamma^- = b > 0$ and find $s_1 < s_0$ and $s_2 > s_0$ so close to $s_0$ that $\sup_{t,u \in [s_1, s_2]} |B_t - B_u| < b/16$. Then $\sup_{x \in \mathbb{Q}} |x - X_{s_2}^{s_1, x}| < b/8$ and $\lambda_{s_2}^- \leq \gamma^- + b/8 < \gamma^+ - b/8 \leq \lambda_{s_2}^+$, by Proposition 2.1. Hence, $\lambda_{s_2}^+ - \lambda_{s_2}^- \geq 3b/4$ and it follows that for some $x \in \mathbb{Q}$, $X_{s_2}^{s_1, x} \in (\lambda_{s_2}^-, \lambda_{s_2}^+)$. This implies that $x \in \Lambda_{s_1}$, so $\Lambda_{s_1} \neq \varnothing$ and $\inf\{t : \Lambda_t \neq \varnothing\} \leq s_1 < s_0$, a contradiction.

Note that we must have $X_1^{s_0,0-} \leq y_3$ and $X_1^{s_0,0+} \geq y_2$. Consider any $x_n \in \Lambda_{s_0+1/n}$. For any fixed $s_3 > s_0$, and $n$ so large that $s_0 + 1/n < s_3$, the processes $\{X_t^{s_0+1/n,x_n}, t \in [s_3,1]\}$ are equicontinuous by Proposition 2.1. This implies that a subsequence converges on $[s_3, 1]$. Since this holds for any $s_3 > s_0$, the diagonal argument can be applied to show that a subsequence of $\{X_t^{s_0+1/n,x_n}, t \in [s_0,1]\}$ converges to some function $\{Y_t, t \in [s_0,1]\}$. One can show that this process solves (1.3) just like in the proof of Proposition 1.1(iii). The three solutions of (1.3) starting from $(s_0, 0)$, $X_t^{s_0,0-}$, $Y_t$ and $X_t^{s_0,0+}$ are distinct because $X_1^{s_0,0-} \leq y_3 < y_1 = Y_1 < y_2 \leq X_1^{s_0,0+}$. □

**3. The lens law $Q$ and anticipated bifurcation times.** We will prove various parts of Theorem 1.6 in different order than stated in the theorem. We will start with the construction and analysis of the "lens law" $Q$. Our first result is a set of explicit formulas we will need in our arguments.

We recall the definition of $Q^{x,y}$. Suppose $x, y \in \mathbb{R}$, $x < y$. Then $Q^{x,y}$ denotes the distribution of $\{(X_t^{0,x}, X_t^{0,y}), t \geq 0\}$ killed at the time $\zeta = \inf\{t > 0 : X_t^{0,x} = X_t^{0,y}\}$. Let $\widehat{L}_t = L_t^{0,x} + L_t^{0,y}$, $\sigma_t = \inf\{s : \widehat{L}_s > t\}$, and $Z_t = |x + \beta L_{\sigma_t}^{0,x} - y - \beta L_{\sigma_t}^{0,y}|$. In other words, $Z_t$ is the distance between $X_t^{0,x}$ and $X_t^{0,y}$ on the time scale defined by the local time clock. Let $\ell = \inf\{t \geq 0 : X_t^{0,x} = X_t^{0,y}\}$ and $\ell_Z = \inf\{t \geq 0 : Z_t = 0\}$. Note that the initial values of $X_t^{0,x}$ and $X_t^{0,y}$ (i.e., $x$ and $y$) are not reflected in the notation for $Z_t$—this is because we will be mostly concerned with the transition probabilities of $Z_t$. Recall that $Z_t$ is a piecewise linear function with the slope $\beta$ or $-\beta$ almost everywhere, on the interval $[0, \ell_Z]$. We have defined $J_t$ to be equal to 0, if at time $\sigma_t$, $X_t^{0,x}$ is at 0, and $J_t = 1$, if at time $\sigma_t$, $X_t^{0,y}$ is at 0. If $\beta > 0$, then $J_t$ is the indicator function of the intervals where $Z_t$ is increasing.

PROPOSITION 3.1. (i) *The process $(Z_t, J_t)$ is Markov with the generator*

$$Af(z,1) = \beta \frac{\partial}{\partial z} f(z,1) - \frac{1+\beta}{2z} f(z,1) + \frac{1+\beta}{2z} f(z,0),$$

$$Af(z,0) = -\beta \frac{\partial}{\partial z} f(z,0) - \frac{1-\beta}{2z} f(z,0) + \frac{1-\beta}{2z} f(z,1).$$

(ii) *The potential density $u((x,j),(z,k))$ of $(Z_t, J_t)$ (i.e., the density of the expectation of the occupation measure) is given by*

$$(3.1) \qquad u((x,0),(z,0)) = \begin{cases} \dfrac{1+\beta}{2\beta^2}, & z \leq x, \\ \dfrac{1-\beta}{2\beta^2} \dfrac{x}{z}, & z > x, \end{cases}$$

$$
(3.2) \qquad u((x,0),(z,1)) = \begin{cases} \dfrac{1-\beta}{2\beta^2}, & z \leq x, \\ \dfrac{1-\beta}{2\beta^2}\dfrac{x}{z}, & z > x, \end{cases}
$$

$$
(3.3) \qquad u((x,1),(z,0)) = \begin{cases} \dfrac{1+\beta}{2\beta^2}, & z \leq x, \\ \dfrac{1+\beta}{2\beta^2}\dfrac{x}{z}, & z > x, \end{cases}
$$

$$
(3.4) \qquad u((x,1),(z,1)) = \begin{cases} \dfrac{1-\beta}{2\beta^2}, & z \leq x, \\ \dfrac{1+\beta}{2\beta^2}\dfrac{x}{z}, & z > x. \end{cases}
$$

PROOF. (i) The claim that $(Z_t, J_t)$ is a Markov process follows easily from the classical Itô excursion theory so we will only outline the evolution of this process and its relationship with $\{(X_t^{0,x}, X_t^{0,y},), t \geq 0\}$. Suppose $J_0 = 0$. As long as $J_t = 0$, only $X_t^{0,x}$ visits 0. The excursions of $X_t^{0,x}$ below 0 occur at the rate $(1-\beta)/2$ and the excursions above 0 occur at the rate $(1+\beta)/2$, on the local time scale, that is, on the same time scale as for the process $Z_t$. When an excursion on the negative side occurs, with the absolute height exceeding $Z_t$, then $J_t$ jumps from 0 to 1 and the analogous process starts: as long as $J_t = 1$, only $X_t^{0,y}$ visits 0. The excursions of $X_t^{0,y}$ below 0 occur at the rate $(1-\beta)/2$ and the excursions above 0 occur at the rate $(1+\beta)/2$, on the local time scale. When an excursion on the positive side occurs, with the height exceeding $Z_t$, then $J_t$ jumps from 1 to 0. The excursions of $X_t^{0,x}$ and $X_t^{0,y}$ are generated according to the same excursion law as for the standard Brownian motion, only their rates are different.

Suppose that $J_0 = 0$. The "probability" that a Brownian excursion has an absolute height greater than $x$ is $1/x$, according to the usual excursion law. It follows from the above that the arrival time of the first negative excursion of $X_t^{0,x}$, with the absolute height exceeding $Z_t$, has the distribution of the first jump arrival time in the Poisson process with variable intensity $(1-\beta)/(2Z_t)$, on the local time scale. This is the same as the jump rate for $J_t$ from 0 to 1. Similarly, the jump rate for $J_t$ from 1 to 0 is the same as the rate of arrival of the first jump in the Poisson process with variable intensity $(1+\beta)/(2Z_t)$. The formula for the generator follows directly from our description of the $(Z_t, J_t)$ evolution.

(ii) It is elementary to check that the function

$$h(z,1) = z, \qquad h(z,0) = \frac{1-\beta}{1+\beta} z,$$

is harmonic for the semigroup of $(Z_t, J_t)$, using the explicit formulas for the generator given in part (i).

For $v > 0$, let $T_v = \inf\{t : Z_t \geq v\}$. We will show that

$$(3.5) \qquad P_{z,0}(T_v < \infty) = \begin{cases} \dfrac{z}{v}\dfrac{1-\beta}{1+\beta}, & v > z, \\ 1, & v \leq z, \end{cases}$$

$$(3.6) \qquad P_{z,1}(T_v < \infty) = \begin{cases} \dfrac{z}{v}, & v > z, \\ 1, & v \leq z. \end{cases}$$

Since $h$ is positive, it follows that $N_t = h(Z_t, J_t)$ is a positive martingale. By the main result in [1], $Z_t \to 0$, a.s., as $t \to \infty$, so $N_t \to 0$, a.s., as $t \to \infty$. By the optional stopping theorem applied to $N_t$ at time $T_v \wedge t$, we get, for $v > z$,

$$\frac{1-\beta}{1+\beta}z = E_{z,0}(N_{T_v \wedge t}) = vP_{z,0}(T_v \leq t) + E_{z,0}(N_t \mathbf{1}_{\{T_v > t\}}).$$

Note that $0 < N_t \mathbf{1}_{\{T_v > t\}} \leq v\mathbf{1}_{\{T_v > t\}}$ and $N_t \mathbf{1}_{\{T_v > t\}} \to 0$ as $t \to \infty$. It follows from the dominated convergence theorem that for $v > z$, $\frac{1-\beta}{1+\beta}z = vP_{z,0}(T_v < \infty)$. For $v \leq z$, $T_v = 0$, $P_{0,z}$-a.s., and this completes the proof of (3.5). We obtain (3.6) in the same manner from

$$z = E_{z,1}(N_{T_v \wedge t}) = vP_{z,1}(T_v \leq t) + E_{z,1}(N_t \mathbf{1}_{\{T_v > t\}}).$$

Suppose that $(Z_0, J_0) = (x, 0)$ and let $S$ be the time of the first jump of $J_t$ to 1. It follows easily from our description of the evolution of $(Z_t, J_t)$ in part (i) of the proof (see also [1]) that

$$P(Z_S/Z_0 < z) = z^{(1-\beta)/(2\beta)}, \qquad 0 \leq z \leq 1.$$

Hence, the density of $Z_S/Z_0$ is $\frac{1-\beta}{2\beta}z^{(1-3\beta)/(2\beta)}$ for $z \in (0,1)$. By the strong Markov property applied at $S$ and (3.6), the probability that $(Z_0, J_0)$ will return to its starting point, that is, $(x, 0)$, is equal to

$$\int_0^1 \frac{1-\beta}{2\beta} z^{(1-3\beta)/(2\beta)} \cdot z \, dz = \frac{1-\beta}{1+\beta}.$$

Hence, the expected number of returns to $(x, 0)$ (including the starting time) is

$$\left(1 - \frac{1-\beta}{1+\beta}\right)^{-1} = \frac{1+\beta}{2\beta}.$$

Since the slope of $Z_t$ is $\beta$, the density of the expectation of the occupation measure at $(x, 0)$ is equal to $(1+\beta)/(2\beta^2)$, as claimed in (3.1). Note that this number does not depend on $x$, so by the strong Markov property applied at the first hitting time of $z < x$, the same formula holds for all $z \leq x$, as stated in (3.1).

If $z > x$, the probability that $Z_t$ will ever hit $z$ is equal to $(x/z) \cdot (1 - \beta)/(1+\beta)$, by (3.5). This combined with the strong Markov property at the hitting time of $z$ and the first part of (3.1) yields the second part of (3.1).

Formula (3.2) can be obtained from (3.1) by noting that for $z > x$, the number of visits to $(z,1)$ is the same as the number of visits to $(z,0)$, a.s. For $z < x$, the number of visits to $(z,1)$ is one less than the number of visits to $(z,0)$.

The other two formulas can be obtained in a very similar manner so the rest of the proof is left to the reader. □

We will give two constructions of $Q$. The first one, presented as a formal proof of Theorem 1.6(ii), is based on an explicit representation of $Q$ in terms of $h$-processes. This construction is followed by a remark containing the second construction, based on Maisonneuve's ideas [13]. The second construction is shorter and has a more abstract character.

PROOF OF THEOREM 1.6. (i) The first assertion is a special case of Lemma 2.8. It is easy to deduce that $G$ is countable from Lemma 2.3(i).

(ii) Recall that the function $h$ given by $h(z,1) = z$, $h(z,0) = z(1-\beta)/(1+\beta)$, is harmonic for $(Z_t, J_t)$. Let $(Z_t^h, J_t^h)$ be the Doob $h$-transform of $(Z_t, J_t)$. The potential density $u^h$ for $(Z_t^h, J_t^h)$ is given by

$$u^h((x,j),(y,k)) = \frac{h(y,k) u((x,j),(y,k))}{h(x,j)}.$$

By Proposition 3.1(ii), for $y > x$, we have $u^h((x,j),(y,k)) = c(j,k)$, where $c(j,k)$ depends only on $\beta$.

Let $(Z_t^{(0,1)}, J_t^{(0,1)})$ denote the process $(Z_t, J_t)$ killed when $Z_t$ escapes from $(0,1)$, and let $(Z_t^{(0,1),h}, J_t^{(0,1),h})$ stand for the process $(Z_t^{(0,1)}, J_t^{(0,1)})$ transformed by $h$. If the process $(Z_t, J_t)$ starts from $(x,j)$ with $x < 1$ and if $Z_t$ leaves $(0,1)$ through 1 at time $s$, then, necessarily, $J_s = 1$. Recall that $h(1,1) = 1$ and $\limsup_{z \to 0} \sup_j h(z,j) = 0$. This implies that $(Z_t^{(0,1),h}, J_t^{(0,1),h})$ is the same as the process $(Z_t, J_t)$ conditioned by the event that $Z_t$ hits 1 before 0.

Let $u^{(0,1)}$ be the potential density for $(Z_t^{(0,1)}, J_t^{(0,1)})$ and let $u^{(0,1),h}$ have the similar meaning for $(Z_t^{(0,1),h}, J_t^{(0,1),h})$. We obviously have $u^{(0,1)} \leq u$ and so $u^{(0,1),h} \leq u^h$. We have shown that $u^h$ is bounded by a constant, so the same applies to $u^{(0,1),h}$. It follows that

$$\sup_{x \in (0,1)} \sum_{k=0,1} \int_0^1 u^{(0,1),h}((x,j),(y,k))\,dy < \infty.$$

This shows that the point $(0,\cdot)$ is an entrance point for $(Z_t^{(0,1),h}, J_t^{(0,1),h})$ [it makes no sense to specify $J_0^h$ when $Z_t^h$ starts from 0 because it is clear that $J_t^h$ will have infinitely many jumps on every interval $(0, \varepsilon)$, $\varepsilon > 0]$. For any $a > 0$, we can construct in a similar way a process $(Z_t^{(0,a),h}, J_t^{(0,a),h})$ starting from $(0, \cdot)$ and representing $(Z_t, J_t)$ conditioned by the event that $Z_t$ hits $a$ before 0.

On some measurable space define probability measures $\widehat{Q}_k$ with disjoint supports, such that $\widehat{Q}_k$ is the distribution of the process $(Z_t^k, J_t^k)$ with the following properties. For an appropriate random variable $T_k \in (0, \infty)$, the process $\{(Z_t^k, J_t^k), t \in [0, T_k]\}$ has the distribution of $(Z_t^{(0,2^{-k}),h}, J_t^{(0,2^{-k}),h})$ starting from $(0, \cdot)$. For $k \geq 1$, the process $\{(Z_t^k, J_t^k), t \geq T_k\}$ starts from $(2^{-k}, 1)$, has the transition probabilities of $(Z_t, J_t)$ conditioned not to hit $(2^{-k+1}, 1)$, and is independent of $\{(Z_t^k, J_t^k), t \in [0, T_k]\}$. We do not have the conditioning in the case $k = 0$, that is, $\{(Z_t^0, J_t^0), t \geq T_0\}$ starts from $(1, 1)$, has the transition probabilities of $(Z_t, J_t)$, and is independent of $\{(Z_t^0, J_t^0), t \in [0, T_0]\}$. Let $\widehat{Q} = \sum_{k \geq 0} 2^k \widehat{Q}_k$. The weights for $\widehat{Q}_k$'s in the sum have been chosen according to (3.6), so that under $\widehat{Q}$, the process $(Z_t^Q, J_t^Q)$ has the same transition probabilities as $(Z_t, J_t)$ on $(s, \infty)$ for any $s > 0$. We add excursions to $(Z_t^Q, J_t^Q)$, in the same manner as described in part (i) of the proof of Proposition 3.1, to obtain a process $\{(X_t^{Q-}, X_t^{Q+}), t \geq 0\}$ with the transition probabilities $Q^{x,y}$ on every interval $(s, \infty)$, $s > 0$. Its distribution will be denoted $Q$.

Next we will argue that $Q$ is the unique measure ("lens law") with the transition probabilities $Q^{x,y}$. Note that given the distribution of $(Z_t^Q, J_t^Q)$, an argument, as in part (i) of the proof of Proposition 3.1, shows that the distribution of $\{(X_t^{Q-}, X_t^{Q+}), t \geq 0\}$ is uniquely defined because under $Q$, the last process has transition probabilities $Q^{x,y}$ on $(s, \infty)$ for every $s > 0$. Let $T_a^Q = \inf\{t : Z_t^Q \geq a\}$. The distribution of $\{(Z_t^Q, J_t^Q), t \geq T_a^Q\}$ is uniquely defined (up to a multiplicative constant) because under $Q$, the transition probabilities for this process are the same as for $(Z_t, J_t)$. The distribution of $\{(Z_t^Q, J_t^Q), t \in [0, T_a^Q]\}$ under $Q$, conditioned by $\{T_a^Q < \infty\}$, is that of $(Z_t^{(0,a),h}, J_t^{(0,a),h})$ starting from $(0, \cdot)$. This concludes the proof of uniqueness for $Q$.

It is clear from the construction of $\widehat{Q}$ and (3.5) that $\lim_{y \downarrow 0}(1/y) Q^{-y,0} = cQ$ for some constant $c$. Recall that we have chosen our normalization of $Q$ so that

$$(3.7) \qquad \lim_{y \downarrow 0}(1/y) Q^{-y,0} = Q.$$

In view of (3.5) and (3.6),

$$(3.8) \qquad \lim_{y \downarrow 0}(1/y) Q^{0,y} = \frac{1-\beta}{1+\beta} Q.$$

(iii) We will only outline the proof of (1.9) because the formula does not require any independence of the lenses and so it is not very deep. Assume without loss of generality that $\beta > 0$. It is enough to prove the formula for sets $A$ of the form $(a, b)$, with $b < 0$ or $a > 0$. Suppose that $a > 0$ and fix some $\varepsilon > 0$. Let $T_1$ be the infimum of times $t$ such that for some lens $U_{s_1}$ with $s_1 \in G$ and $B_{s_1} \in (a, b)$, we have $L_t^{s_1, 0+} \geq L_t^{s_1, 0-} + \varepsilon$. Note that $T_1$ is a stopping time. Clearly, $\{U_{s_1}(t), t \geq T_1 - s_1\}$ has the same transition probabilities as $Q$ and its value at time $T_1 - s_1$ is $(-\varepsilon, 0)$. Similarly, let $T_k$ be the infimum of times $t > T_{k-1}$ such that for some lens $U_{s_k}$ with $s_k \in G$ and $B_{s_1} \in (a, b)$, we have $L_t^{s_k, 0+} \geq L_t^{s_k, 0-} + \varepsilon$. We see that $T_k$ is a stopping time and $\{U_{s_k}(t), t \geq T_k - s_k\}$ has the same transition probabilities as $Q$. Summing over all $k$ and using (3.7), we obtain (1.9) for functions $f$ which depend only on the post-$T$ process, where $T$ is the infimum of times such that $Z_t = \varepsilon$. The general result is obtained by letting $\varepsilon \to 0$.

The case $b < 0$ requires an application of (3.8) instead of (3.7) because if we follow an analogous argument, we have $U_{s_1}(T_1 - s_1) = (0, \varepsilon)$ and not $U_{s_1}(T_1 - s_1) = (-\varepsilon, 0)$. This explains the factor $\frac{1-\beta}{1+\beta}$ on the left-hand side of (1.9). □

REMARK 3.1. We will now sketch an alternative construction of the lens law $Q$. We need the usual general Markov process setup, with some probability space $(\Omega, \mathcal{F}, P)$, filtration $\{\mathcal{F}_t\}$ and shift operators $\theta_t$ on $\Omega$ that act on the Brownian motion $B_t$ as usual, that is, $B_s(\theta_t) = B_{s+t} - B_t$, and on $X_t^x$ in the following manner: $X_s^x(\theta_t) = X_s^{X_t^x}(\theta_t)$. Consider the increasing process

$$A_t = \sum_{\substack{s \in G \\ s \leq t}} \int_0^{L(\theta_s)} e^{-u} \tfrac{1}{2}(e^{-|\bar{\ell}(u)|} - e^{-|\bar{u}(u)|}) \, du,$$

where $L$ is the length (in time units) of the lens, $\bar{\ell}(u)$ is its lower limit at time $u$ and $\bar{u}(u)$ is its upper limit. Since $\int_{t=0}^{\infty} e^{-t} \tfrac{1}{2} \int_{-\infty}^{\infty} e^{-|x|} \, dx \, dt = 1$, and in the half-plane $\mathbb{R}_+ \times \mathbb{R}$, lenses do not intersect each other (although they may touch), we see that $\int_0^{\infty} e^{-t} \, dA_t \leq 1$. It can be easily verified that $A_{t+s} = A_t + A_s \circ \theta_t$. Thus, $(A_t)$ is a raw additive functional of $\mathcal{F}_t$, and since no $(\mathcal{F}_t)$ stopping time passes through its jumps, its dual predictable projection $(\tilde{A}_t)$ is a continuous additive functional of $(\mathcal{F}_t)$, which we shall call the lens local time. The functional $x \to X_\cdot^x$ is monotone, in the usual sense of inequality between functions. The family of such functionals is "good" so arguing as in [13], we can prove that there exists a kernel $Q^{X_\cdot}$ such that, for every $(\mathcal{F}_t)$-predictable $V$ and $\mathcal{F}$-measurable $f$,

$$E \sum_{s \in G} V_s f \circ \theta_s = E \int_0^{\infty} V_s Q^{X_s(\cdot)}(f) \, d\tilde{A}_s.$$

Note further that the evolution of a lens (its width, the position of its upper and lower parts during times between its formation and its coalescence, etc.)

is independent of the value of the function $X_s^{\cdot}$ at time $s$ when the lens starts. It only depends on the Brownian motion $(B_t \circ \theta_s)_{t \geq 0}$ and the fact that at time $s$, $X_s^x = 0$ for some $x$, and for a sequence $x_n$ increasing in $n$, with $X_s^{x_n} < 0$ for all $n$, $X_s^{x_n}$ increases to 0 as $n \to \infty$. All points $s \in G$ satisfy this condition and, therefore, $\tilde{A}$ does not charge points for which $X_s^{\cdot}$ does not satisfy this condition. In particular, $\int_0^{L(\theta_s)} e^{-u} \frac{1}{2}(e^{-|\bar{\ell}(u)|} - e^{-\bar{u}(u)})\, du$ depends only on $(B_t \circ \theta_s)_{t \geq 0}$. It follows, as in the derivation of "exit systems" in [13], that if $f$ is a function of the lens only and is independent of the particular form of $X^{\cdot}(\theta_s)$, then $Q^{X_s}(f)$ does not depend on $X_s^{\cdot}$ and we denote it $Q(f)$. The normalization given in Theorem 1.6(ii) corresponds to the following normalization of $\tilde{A}$ and $Q$ in the present context,

$$\tilde{A}_t = \frac{1-\beta}{1+\beta}\left(-\min_{u \leq t} B_u\right) + \max_{u \leq t} B_u$$

and $Q(l > v) = 1/v$, where $l$ is the maximal opening of the lens.

The Williams decomposition of the Brownian excursion can be presented as follows. Fix some $b > 0$. Let $R_t^1$ and $R_t^2$ be two independent Brownian motions starting from 0 and conditioned to go to infinity before returning to 0 (i.e., they are three-dimensional Bessel processes). Kill $R_t^1$ and $R_t^2$ when they hit $b$, call the resulting processes $\widetilde{R}_t^1$ and $\widetilde{R}_t^2$, time-reverse $\widetilde{R}_t^2$ to obtain $\widehat{R}_t^2$, and concatenate $\widetilde{R}_t^1$ and $\widehat{R}_t^2$. The result is a process with the same distribution as the Brownian excursion conditioned to have height $b$. A very similar construction of $Z_t$ under $Q$, conditioned to have height $b$, can be given. The following is an informal version of the argument given in the next proof, with minor changes. Let $V_t^1$ and $V_t^2$ be two independent processes with the same transition probabilities as those of $Z_t$ under $Q^{x,y}$, except that $|V_t^1 - V_t^2|$ is conditioned to go to infinity before returning to 0, and assume that $V_t^1$ and $V_t^2$ start from 0. Kill $V_t^1$ and $V_t^2$ when $|V_t^1 - V_t^2| = b$, call the resulting processes $\widetilde{V}_t^1$ and $\widetilde{V}_t^2$, time-reverse $\widetilde{V}_t^2$ to obtain $\widehat{V}_t^2$, and concatenate $\widetilde{V}_t^1$ and $\widehat{V}_t^2$. This construction yields a process with the same distribution as $Z_t$ under $Q$, conditioned to have height $b$.

PROOF OF THEOREM 1.6. (vi) We will use Nagasawa's theorem on time reversal (see [14] or [17]). In order to do that, we need to calculate generators and potentials for some processes and find a reference measure under which the processes are dual.

Fix any $b > 0$, let $T_b = \inf\{t : Z_t = b\}$, and recall the generator $A$ and a harmonic function $h$ from the statement of Proposition 3.1 and part (ii) of its proof. Let $(Z_t^{(0,b)}, J_t^{(0,b)})$ denote the process $(Z_t, J_t)$ killed when $Z_t$ exits $(0,b)$, and let $(Z_t^{(0,b),h}, J_t^{(0,b),h})$ be $(Z_t^{(0,b)}, J_t^{(0,b)})$ conditioned by $h$. The process $(Z_t^{(0,b),h}, J_t^{(0,b),h})$ is $(Z_t, J_t)$ conditioned by $\{T_b < \infty\}$.

Let $h^-(x,1) = 1 - x/b$, $h^-(x,0) = 1 - (x/b)(1-\beta)/(1+\beta)$ and note that $h^-$ is harmonic for $A$, with the boundary value 0 at $(b,1)$ and 1 at $(0,j)$. Let $(Z_t^{(0,b),h^-}, J_t^{(0,b),h^-})$ be $(Z_t^{(0,b)}, J_t^{(0,b)})$ conditioned by $h^-$ and note that $(Z_t^{(0,b),h^-}, J_t^{(0,b),h^-})$ is $(Z_t, J_t)$ conditioned by $\{T_b = \infty\}$.

Consider $(Z_t, J_t)$ under $Q$ conditioned by $\{\sup_{t \geq 0} Z_t = b\}$ and note that the distribution of $\{(Z_t, J_t), t \in [0, T_b]\}$ under $Q$ is the same as the distribution of $(Z_t^{(0,b),h}, J_t^{(0,b),h})$ starting from $(0, \cdot)$, while the distribution of $\{(Z_t, J_t), t \in [T_b, T_0]\}$ is that of the process $(Z_t^{(0,b),h^-}, J_t^{(0,b),h^-})$ starting from $(b,0)$. It will suffice to show that the time reversal of $Z_t^{(0,b),h^-}$ starting from $b$ is the same, in distribution, as $Z_t^{(0,b),h}$ starting from 0.

By Proposition 3.1(i), and using conditioning by $h$, the generator of $(Z_t^{(0,b),h}, J_t^{(0,b),h})$ is equal to

$$(3.9) \quad A_\beta^h f(z,1) = \beta \frac{\partial}{\partial z} f(z,1) - \frac{1-\beta}{2z} f(z,1) + \frac{1-\beta}{2z} f(z,0),$$

$$(3.10) \quad A_\beta^h f(z,0) = -\beta \frac{\partial}{\partial z} f(z,0) - \frac{1+\beta}{2z} f(z,0) + \frac{1+\beta}{2z} f(z,1).$$

Consider a process $(\hat{Z}_t, \hat{J}_t)$ defined in the same manner as $(Z_t, J_t)$, except that the skewness parameter should be $-\beta$ instead of $\beta$. Recall that $J_t$ is the indicator function of the intervals where $Z_t$ is increasing. The process $\hat{J}_t$ is the indicator function of the intervals where $\hat{Z}_t$ is decreasing. We write $A_\beta$ instead of $A$ to emphasize the dependence of the generator $A$, defined in Proposition 3.1, on the parameter $\beta$. By that result, $\hat{A} = A_{-\beta}$. The function $\hat{h}(z,0) = z$, $\hat{h}(z,1) = \frac{1-\beta}{1+\beta} z$ is harmonic for $(\hat{Z}_t, \hat{J}_t)$, and the generator $\hat{A}^{\hat{h}}$ of the $\hat{h}$-transform $(\hat{Z}_t^{\hat{h}}, \hat{J}_t^{\hat{h}})$ of $(\hat{Z}_t, \hat{J}_t)$ satisfies $\hat{A}^{\hat{h}} = A_{-\beta}^h$. This differs from (3.9) and (3.10) only in that the roles of the states 0 and 1 of $J_t$ have been reversed. This means that the distribution of $Z_t^{(0,b),h}$ starting from 0 and the distribution of $\hat{Z}_t^{\hat{h}}$ starting from 0 and killed at $b$ are identical. It remains to show that the time reversal of $Z_t^{(0,b),h^-}$ starting from $b$ has the same distribution as $\hat{Z}_t^{\hat{h}}$ starting from 0 and killed at $b$.

In view of Proposition 3.1(i), the generator of $(Z_t^{(0,b),h^-}, J_t^{(0,b),h^-})$ is equal to

$$A^{h^-} f(z,1) = \beta \frac{\partial}{\partial z} f(z,1) - \frac{1}{2z}\left(1 + \beta \frac{b+z}{b-z}\right) f(z,1)$$
$$+ \frac{1}{2z}\left(1 + \beta \frac{b+z}{b-z}\right) f(z,0),$$

$$A^{h^-} f(z,0) = -\beta \frac{\partial}{\partial z} f(z,0) - \frac{1-\beta^2}{2z} \frac{b-z}{b-z+\beta(b+z)} f(z,0)$$

$$+ \frac{1-\beta^2}{2z} \frac{b-z}{b-z+\beta(b+z)} f(z,1).$$

Let $m$ be defined by

$$m(dz,1) = \frac{1-\beta}{2\beta^2}\left(1-\frac{z}{b}\right) dz,$$

$$m(dz,0) = \left(\frac{1+\beta}{2\beta^2} - \frac{1-\beta}{2\beta^2}\frac{z}{b}\right) dz,$$

and $(f,g)_m = \int fg\, dm$. Note that $m$ depends on $b$. We want to show that $\hat{A}^{\hat{h}}$ and $A^{h^-}$ are in duality with respect to $m$, that is,

(3.11) $$(\hat{A}^{\hat{h}} g, f)_m = (g, A^{h^-} f)_m,$$

for all $C^1$-functions $f,g$ with compact support in $(0,b)$. We omit tedious but completely elementary calculations which show that the last formula is, indeed, true.

Next we will show that the duality measure $m$ is the potential of $(\hat{Z}_t^{\hat{h}}, \hat{J}_t^{\hat{h}})$ starting at $(0,0)$ (or rather when $\hat{Z}_0^{\hat{h}} = 0$) and killed when $\hat{Z}_t^{\hat{h}}$ exits $(0,b)$.

Recall the function $u$ defined in Proposition 3.1(ii). We will make its dependence on $\beta$ explicit by writing $u_\beta$. When we apply Proposition 3.1(ii) with $-\beta$ in place of $\beta$, we see that the potential density $\hat{u}((x,j),(z,k))$ of $(\hat{Z}_t, \hat{J}_t)$ is given by $\hat{u} = u_{-\beta}$. Since $\hat{h}(z,0) = z$, $\hat{h}(z,1) = \frac{1-\beta}{1+\beta} z$, this implies that the potential density $\hat{u}^{\hat{h}}((x,j),(z,k))$ of $(\hat{Z}_t^{\hat{h}}, \hat{J}_t^{\hat{h}})$ is equal to

$$\hat{u}^{\hat{h}}((x,0),(z,0)) = \begin{cases} \dfrac{1-\beta}{2\beta^2}\dfrac{z}{x}, & z \leq x, \\ \dfrac{1+\beta}{2\beta^2}, & z > x, \end{cases}$$

$$\hat{u}^{\hat{h}}((x,0),(z,1)) = \begin{cases} \dfrac{1-\beta}{2\beta^2}\dfrac{z}{x}, & z \leq x, \\ \dfrac{1-\beta}{2\beta^2}, & z > x, \end{cases}$$

$$\hat{u}^{\hat{h}}((x,1),(z,0)) = \begin{cases} \dfrac{1+\beta}{2\beta^2}\dfrac{z}{x}, & z \leq x, \\ \dfrac{1+\beta}{2\beta^2}, & z > x, \end{cases}$$

$$\hat{u}^{\hat{h}}((x,1),(z,1)) = \begin{cases} \dfrac{1+\beta}{2\beta^2}\dfrac{z}{x}, & z \leq x, \\ \dfrac{1-\beta}{2\beta^2}, & z > x. \end{cases}$$

Let $(\hat{Z}_t^{(0,b),\hat{h}}, \hat{J}_t^{(0,b),\hat{h}})$ denote the process $(\hat{Z}_t^{\hat{h}}, \hat{J}_t^{\hat{h}})$ killed when $\hat{Z}_t^{\hat{h}}$ exits $(0,b)$, and let $\hat{u}^{(0,b),\hat{h}}((x,j),(z,k))$ denote the corresponding potential density. For $x < b$, the process $(\hat{Z}_t^{\hat{h}}, \hat{J}_t^{\hat{h}})$ hits $(b,1)$ with probability 1, so by the strong Markov property applied at the hitting time of $(b,1)$, for $x, z \in (0,b)$,

$$\hat{u}^{(0,b),\hat{h}}((x,j),(z,k)) = \hat{u}^{\hat{h}}((x,j),(z,k)) - \hat{u}^{\hat{h}}((b,1),(z,k)).$$

It is now straightforward to verify that

$$m(dz, 1) = \frac{1-\beta}{2\beta^2}\left(1 - \frac{z}{b}\right) dz = \hat{u}^{(0,b),\hat{h}}((0,0),(z,1))\, dz,$$

$$m(dz, 0) = \left(\frac{1+\beta}{2\beta^2} - \frac{1-\beta}{2\beta^2}\frac{z}{b}\right) dz = \hat{u}^{(0,b),\hat{h}}((0,0),(z,0))\, dz.$$

This completes the proof that the duality measure $m$ is the potential of $(\hat{Z}_t^{(0,b),\hat{h}}, \hat{J}_t^{(0,b),\hat{h}})$ starting at $(0,0)$. We have already shown that $\hat{A}^{\hat{h}}$ and $A^{h^-}$ are dual [see (3.11)], so Nagasawa's theorem implies that the time reversal of $Z_t^{(0,b),h^-}$ starting from $b$ has the same distribution as $\hat{Z}_t^{(0,b),\hat{h}}$ starting from 0. We have already shown that the laws of $Z_t^{(0,b),h}$ and of $\hat{Z}_t^{(0,b),\hat{h}}$ are the same. This completes the proof of part (vi) of Theorem 1.6.

(vii) Recall the argument in part (i) of the proof of Proposition 3.1 and the notation in the paragraph following the statement of Theorem 1.6. Suppose that $Z_t$ is increasing when $X_t$ is hitting 0 (this depends on the sign of $\beta$). After time reversal, $Z_t$ is decreasing when (the time-reversal of) $X_t$ is hitting 0. Since $Z_{\widehat{L}_t}$ is adapted to the filtration generated by $(X_t, \widetilde{X}_t)$, we see that $Q$-lenses are not invariant under time reversal.

(v) We will discuss only the case $\beta > 0$. For $x < y$, let

$$V(x,y) = \inf\{z > y : x + \beta L_t^x = y + \beta L_t^y = z \text{ for some } t\}$$

and let $Z_t^{x,y}$ be defined relative to $X_t^x$ and $X_t^y$ in the same way as $Z_t^s$ in Theorem 1.6. It has been shown in the proof of Theorem 1.1 in [4] that if $x < y < z$, then $V(x,y)$ is independent of $V(y,z)$. The proof is based on the following idea. The value of $V(x,y)$ is determined by the excursions of $X_t^y$ above 0, and the value of $V(y,z)$ is determined by excursions of $X_t^y$ below zero. The two excursion processes are independent, so $V(x,y)$ and $V(y,z)$ are independent. The same argument shows that the processes $Z_{\cdot}^{x,y}$ and $Z_{\cdot}^{y,z}$ are independent, for $x < y < z$, where $Z_{\cdot}^{x,y}$ is the process defined at the beginning of the proof. We can generalize this as follows. If $x_1 < x_2 < \cdots < x_n$, then the processes $\{Z_{\cdot}^{x_k, x_{k+1}}\}_{1 \leq k \leq n-1}$ are jointly independent.

For any integer $n$, let $\mathcal{A}_n$ be the collection of pairs $(k2^{-n}, \{Z_{\cdot}^{k2^{-n},(k+1)2^{-n}}\})_{k \in \mathbb{Z}}$. The independence of $Z_{\cdot}^{k2^{-n},(k+1)2^{-n}}$'s and (3.7) and (3.8) imply easily that the point processes $\mathcal{A}_n$ converge weakly to a Poisson point process on the space $\mathbb{R} \times C[0,\infty)$ and intensity

$$(3.12) \qquad \mathbf{1}_{(-\infty,0]}(x)\,dx \times \frac{1-\beta}{1+\beta}Q^Z + \mathbf{1}_{(0,\infty)}(x)\,dx \times Q^Z,$$

where $Q^Z$ is the $Q$-distribution of $Z$.

Let $C_\ell[0,\infty)$ be the set of functions $f \in C[0,\infty)$ for which $\inf\{t > 0 : f(t) = 0\} \geq \ell$. We will show that $\mathcal{A}_n$'s converge almost surely on $\mathbb{R} \times C_\ell[0,\infty)$ for every $\ell > 0$. Let $\mathcal{A}_n(\ell)$ be the number of points of $\mathcal{A}_n$ in $[0,1) \times C_\ell[0,\infty)$. It was noticed in [4] that for $x < y < z$, we have $\max(V(x,y), V(y,z)) = V(x,z)$. This implies that $\mathcal{A}_{n+1}(\ell - 2^{-n-1}) \geq \mathcal{A}_n(\ell)$ and so the limit $\mathcal{A}(\ell) = \lim_{n\to\infty} \mathcal{A}_n(\ell + 2^{-n})$ is well defined, a.s. We also have another monotonicity property, namely, $Z_t^{2k2^{-n-1},(2k+1)2^{-n-1}} \leq Z_t^{k2^{-n},(k+1)2^{-n}}$ for all $t$ and, similarly, $Z_t^{(2k+1)2^{-n-1},(2k+2)2^{-n-1}} \leq Z_t^{k2^{-n},(k+1)2^{-n}}$. It is not hard to deduce from this that for some point process $\mathcal{A}$, the processes $\mathcal{A}_n$'s converge almost surely to $\mathcal{A}$ on $[0,1) \times C_\ell[0,\infty)$. Clearly, $\mathcal{A}$ is a Poisson point process whose distribution is given by (3.12) restricted to $[0,1) \times C_\ell[0,\infty)$. An extension to $[0,\infty) \times C[0,\infty)$ is routine. The extra factor in the formula for the intensity on $(-\infty,0] \times C[0,\infty)$ can be justified as in the proof of Theorem 1.6(iii) [see (3.7) and (3.8)].

Recall that the equation (1.3) has unique strong solutions for all rational $s$ and $x$ simultaneously, a.s. Hence, there are no anticipated bifurcation times $s$ with $B_s$ rational. If $x$ is irrational then for every $n > 1$, it belongs to an interval $(k2^{-n}, (k+1)2^{-n})$. It follows easily from the definition of a bifurcation time and that of $Z_t^{k2^{-n},(k+1)2^{-n}}$ that there exists an anticipated bifurcation time $s$ with $B_s = -x$ if and only if $(x,f) \in \mathcal{A}$ for some $f \in C_\ell[0,\infty)$ with $\ell > 0$.

(iv) This part of the theorem is not much more than a "soft" remark so we will only sketch the proof. The idea of the argument is that different lenses corresponding to anticipated bifurcation times overlap on the time scale and so they cannot be independent because Brownian paths generated independently cannot have the same shape.

Recall that $B_t$ denotes a one-dimensional Brownian motion. The three-dimensional Brownian motion does not hit a fixed straight line, a.s. This and standard arguments easily imply that for any fixed $s_1, s_2 > 0$ and $x_1, x_2 \in \mathbb{R}$, there are no $t_1, t_2, t_3 \geq 0$ such that $t_2 - t_1 = s_1$, $t_3 - t_2 = s_2$, $B_{t_2} - B_{t_1} = x_1$, and $B_{t_3} - B_{t_2} = x_2$. Let $B'_t$ be a Brownian motion independent of $B_t$. It follows by conditioning on the values of $B'_{t_2} - B'_{t_1}$ and $B'_{t_3} - B'_{t_2}$ that if $t_1, t_2, t_3 \geq 0$ are fixed, then with probability 1, there is no $u$ such that $B'_{t_2} - B'_{t_1} = B_{t_2+u} - B_{t_1+u}$ and $B'_{t_3} - B'_{t_2} = B_{t_3+u} - B_{t_2+u}$. This holds for all rational triplets $t_1, t_2, t_3 \geq 0$ simultaneously, a.s.

Suppose that $\mathcal{B}$ is a Poisson process, that is, anticipated lenses are generated independently. It is easy to see that the upper part of some anticipated lens agrees with the lower part of some other anticipated lens on a nonempty open interval. Since such an interval contains a nonempty subinterval on which neither process visits 0, we see that two Brownian paths generated independently agree on a nondegenerate interval. This contradiction shows that $\mathcal{B}$ is not Poisson. $\square$

PROOF OF THEOREM 1.5(ii). Fix an integer $n > 1$, let $\varepsilon = 1/n$ and $s_k = k/n$ for $k = 0, 1, \ldots, n$. Note that if four different solutions to (1.3) start at $(s, x)$, then necessarily $x = 0$. Fix $\alpha \in (2, 3)$. Suppose without loss of generality that $\beta > 0$ and for a lens $\{(X_t^{s,0-}, X_t^{s,0+}), t \in [s, u]\}$, let $\chi(s) = \beta L_u^{s,0+} = \beta L_u^{s,0-}$. Suppose that four distinct solutions start at a point $(s, 0)$. We will write $\chi_1(s)$ for the quantity analogous to $\chi(s)$, but defined relative to the lowest two solutions starting at $(s, 0)$. Similarly, $\chi_2(s)$ will correspond to the middle two solutions, and $\chi_3(s)$ to the upper two. Fix some integer $k$. Let $A_k$ be the event that four different solutions start at $(s, 0)$ with $s \in [s_k, s_{k+1})$ and $\chi_j(s) > 1$ for $j = 1, 2, 3$. We will partition $A_k$ into two events, $A_k^1$ and $A_k^2$. The event $A_k^1$ is when $A_k$ occurs and $\sup_{t \in [s_{k-1}, s_k]} |B_t - B_{s_{k-1}}| \geq \varepsilon^{1/\alpha}$. By (2.8),

$$(3.13) \qquad P(A_k^1) \leq c_3 \varepsilon^{1/\alpha - 1/2} \exp(-\tfrac{1}{2}\varepsilon^{2/\alpha - 1}).$$

We let $A_k^2 = A_k \setminus A_k^1$. Let $x_1, x_2, x_3$ and $x_4$ be the values of the four solutions starting at $(s, 0)$, at the time $s_{k+1}$. Since the points $(s_{k+1}, x_j), j = 1, 2, 3, 4$, lie on some solutions to (1.3), away from their starting points, they actually lie on some solutions to (1.3) with rational coordinates of starting points. Hence, $X_\cdot^{s_{k+1}, x_j-} \equiv X_\cdot^{s_{k+1}, x_j+}$ and we will denote these processes $X_\cdot^{s_{k+1}, x_j}$. If $A_k^2$ holds, then $|x_k| \leq 2\varepsilon^{1/\alpha}$ for $k = 1, 2, 3, 4$, by the definition of $A_k^1$ and Proposition 2.1. Recall from the proof of Proposition 2.1 that $\{\beta L_u^{t_1, x}, t_1 \leq u \leq t_2\}$ has the same modulus of continuity as $\{B_u, t_1 \leq u \leq t_2\}$, for any $t_1 < t_2$. Assume that $\varepsilon$ is so small that $\varepsilon^{1/\alpha} < \tfrac{1}{4}$. Since $\chi_j(s) > 1$ for every $j$, it follows that $x_j + \beta L_\cdot^{s_{k+1}, x_j}$ must increase by at least $\tfrac{1}{2}$ before it meets any other $x_m + \beta L_\cdot^{s_{k+1}, x_m}, m \neq j$. Hence, there exist (at least) three anticipated bifurcation times with $\chi > \tfrac{1}{2}$, relative to time $s_{k+1}$ rather than time 0, that is, for some $y_1 < y_2 < y_3$, and every $j = 1, 2, 3$, we have $\inf\{z : z = y_j + \beta L_t^{s_{k+1}, y_j-} = y_j + \beta L_t^{s_{k+1}, y_j+}$ for some $t\} \geq y_j + \tfrac{1}{2}$. Theorem 1.6 can be applied to processes starting from $s_{k+1}$ by shift invariance of the Brownian motion. Note that $|y_j| \leq 2\varepsilon^{1/\alpha}$, because a similar bound holds for $x_j$'s. We see that if $A_k^2$ occurs, then there are three anticipated bifurcation times corresponding to solutions of (1.3), starting at time $s_{k+1}$, at space points within $2\varepsilon^{1/\alpha}$ of 0. It is easy to see that $Q(\chi > 1) < \infty$. Note that $\chi$ is a function of $Z$, so by part (v) of Theorem 1.6, the "$\chi$" point process is Poisson and, therefore, the probability of $A_k^2$ is bounded by $c_1(2\varepsilon^{1/\alpha})^3$. Hence, the probability of $\bigcup_{0\leq k\leq n-1} A_k^2$ is bounded by $c_2(2\varepsilon^{1/\alpha})^3\varepsilon^{-1} = c_3\varepsilon^{3/\alpha-1}$. Since $\alpha \in (2,3)$, this goes to zero as $\varepsilon \to 0$. This combined with (3.13) yields $P(\bigcup_{0\leq k\leq n-1} A_k) \to 0$ as $n = 1/\varepsilon \to \infty$. We conclude that there are no $s \in [0,1)$, where four different solutions start with $\chi_j(s) > 1$ for every $j$. The same argument shows that for every integer $m \geq 1$, there are no $s \in [-m,m)$, where four different solutions start with $\chi_j(s) > 1/m$ for every $j$, a.s. Since $m$ is arbitrarily large, the proof is complete. □

**Acknowledgments.** We would like to thank Professor Boris Tsirelson for his instructive pictures and the very helpful discussions that came along with them. The first author is grateful for the opportunity to visit Technion (Haifa) where part of this paper was written. We are grateful to the referees for a very careful reading of the original manuscript and many useful suggestions.

DEPARTMENT OF MATHEMATICS
UNIVERSITY OF WASHINGTON
SEATTLE, WASHINGTON 98195-4350
USA
E-MAIL: burdzy@math.washington.edu
URL: www.math.washington.edu/~burdzy/

FACULTY OF INDUSTRIAL ENGINEERING
 AND MANAGEMENT
TECHNION INSTITUTE
HAIFA 32000
ISRAEL
E-MAIL: iehaya@tx.technion.ac.il
URL: iew3.technion.ac.il/Home/users/iehaya.phtml?YF